\documentclass[11pt]{article}
\usepackage{amscd,amsmath,amssymb,amsthm,latexsym}

\input rlepsf

\def\co{\colon\thinspace}
\newcommand{\SL}{\mbox{\rm SL}}
\newcommand{\GL}{\mbox{\rm GL}}
\newcommand{\E}{\mbox{\rm E}}
\newcommand{\PSL}{\mbox{\rm PSL}}
\newcommand{\SO}{\mbox{\rm SO}}
\newcommand{\ST}{\mbox{\rm ST}}

\newcommand{\SU}{\mbox{\rm SU}}
\newcommand{\Diff}{\mbox{\rm Diff}}
\newcommand{\Inn}{\mbox{\rm Inn}}
\newcommand{\Aut}{\mbox{\rm Aut}}
\newcommand{\Out}{\mbox{\rm Out}}
\newcommand{\Hom}{\mbox{\rm Hom}}
\newcommand{\Isom}{\mbox{\rm Isom}}
\newcommand{\Spin}{\mbox{\scriptsize\rm Spin}}
\newcommand{\T}{{\mathcal T}}
\newcommand{\G}{{\mathcal G}}
\newcommand{\M}{{\mathcal M}}
\newcommand{\R}{{\mathcal R}}
\newcommand{\C}{{\mathcal C}}
\newcommand{\B}{{\mathcal B}}
\newcommand{\F}{{\mathcal F}}
\newcommand{\alg}{\mbox{\scriptsize\rm alg}}
\newcommand{\orb}{\mbox{\scriptsize\rm orb}}
\newcommand{\ou}{\overline{u}}
\newcommand{\ov}{\overline{v}}
\newcommand{\oq}{\overline{q}}
\newcommand{\orh}{\overline{\rho}}
\newcommand{\oph}{\overline{\varphi}}
\newcommand{\osi}{\overline{\sigma}}

\newtheorem{thm}{Theorem}[section]
\newtheorem{lem}[thm]{Lemma}
\newtheorem{prop}[thm]{Proposition}

\newtheorem{defn}[thm]{Definition}
\newtheorem{rem}[thm]{Remark}

\begin{document}

\title{Moduli of Contact Circles}
\author{Hansj\"org Geiges and Jes\'us Gonzalo}
\date{}

\maketitle

{\footnotesize \tableofcontents}


\vspace{1cm}


This paper\footnote{{\it 2000 MSC:}
53D35 (Primary); 58D27, 32G15, 57S30 (Secondary)}
is a continuation of our study of contact circles
begun in~\cite{gego95}. In that paper we gave a complete
classification of the closed, orientable $3$-manifolds that admit what
we called a {\em taut contact circle}. Below we recall the definition of
this type of structure and some fundamental concepts and results
from our previous paper.

The main objective of the present paper is to introduce and describe
deformation spaces for taut contact circles. In~\cite{gego95} we had
already observed that taut contact circles admit non-trivial deformations,
and we had described the moduli spaces of taut contact circles on spherical
$3$-manifolds. Here we introduce Teichm\"uller and moduli spaces
for taut contact circles in a more systematic fashion, and we completely
determine these spaces for all 3-manifolds admitting taut contact circles.
\section{Definitions and previous results}
In this section we briefly summarise the basic definitions and the main
result from~\cite{gego95}. In later sections we shall occasionally
quote other results from our earlier paper without much further explanation,
but on the whole we have tried to make the present paper reasonably
self-contained.

\begin{defn}
\label{defn:circle}
A {\bf taut contact circle} on a $3$-manifold $M$ is a pair of
contact forms $(\omega_1,\omega_2)$ such that the $1$-form $\lambda_1\omega_1
+\lambda_2\omega_2$ is a contact form defining the same volume form
for all $(\lambda_1,\lambda_2) \in S^1\subset {\mathbb R}^2$. Equivalently,
we require that the following equations be satisfied:
\begin{eqnarray*}
\omega_1\wedge d\omega_1 & = & \omega_2\wedge d\omega_2 \neq 0,\\
\omega_1\wedge d\omega_2 & = & -\omega_2\wedge d\omega_1.
\end{eqnarray*}
If the mixed terms $\omega_1\wedge d\omega_2$ and $\omega_2\wedge d\omega_1$
are identically zero rather than just of opposite sign, we speak of a
{\bf Cartan structure}.
\end{defn}

The following was our main classification result in~\cite{gego95}.

\begin{thm}
\label{thm:classification}
Let $M$ be a closed $3$-manifold. Then $M$ admits
a taut contact circle
if and only if $M$ is diffeomorphic to a quotient of the Lie group
$\G$ under a discrete subgroup $\Gamma$ acting by left multiplication,
where $\G$ is one of the following.

{\em (a)} $S^{3}=\SU (2)$, the universal cover of $\SO(3)$.

{\em (b)} $\widetilde{\SL}_{2}$, the universal cover of $\PSL_{2}{\mathbb R}$.

{\em (c)} $\widetilde{\E}_{2}$, the universal cover of the Euclidean
group (that is, orientation preserving isometries of the
Euclidean plane).

Each of these manifolds admits a Cartan structure.
\end{thm}

If $(\omega_1,\omega_2)$ is a taut contact circle, then so is
$(v\omega_1,v\omega_2)$ for any positive function $v$, and so is
$(\omega_1\cos\theta -\omega_2\sin\theta ,\omega_1\sin\theta +\omega_2
\cos\theta )$ for any constant angle~$\theta$. When classifying
taut contact circles, we shall do so up to this ambiguity. To this end we
introduce the following concepts.

\begin{defn}
The {\bf conformal class} (resp.\ {\bf homothety class})
of a taut contact circle $(\omega_1,\omega_2)$ is the collection of all
pairs $(\omega_1',\omega_2')$ obtained from $(\omega_1,\omega_2)$ by
multiplication by some positive function~$v$ (resp.\ multiplication by
some~$v$ and rotation by some~$\theta$).
\end{defn}
\section{Deformation spaces}
\label{section:deform}
We now introduce Teichm\"uller and moduli spaces for taut contact circles
in complete analogy to the corresponding definitions for complex or
other geometric structures. Let $M$ be a manifold as in
Theorem~\ref{thm:classification}. Write $\Diff (M)$ for the diffeomorphism
group of $M$ and $\Diff_0(M)$ for the group of diffeomorphisms isotopic
to the identity.
Let ${\mathcal C}(M)$ be the space of homothety classes
of taut contact circles on~$M$. Our main objects of study in the
present paper are the following spaces.

\begin{defn}
The {\bf Teichm\"uller space} of taut contact circles on $M$ is
${\T}(M)={\mathcal C}(M)/\Diff_0 (M)$.

The {\bf moduli space} of taut contact circles is ${\M}(M)=
{\mathcal C}(M)/\Diff(M)$.
\end{defn}

In terms of this language, the results from~\cite[Section~6]{gego95}
about the moduli of taut contact circles on left-quotients of~$S^3$
can be stated as follows.

\begin{thm}
\label{thm:deform-s}
Let $M=\Gamma\backslash \SU (2)$ be a left-quotient of $\SU(2)$ under
a (discrete, cocompact) subgroup of $\SU (2)$. If $\Gamma$
is non-abelian, then
${\M}(M)$ consists of a single point. Otherwise, $M$ must be
a lens space $L(m,m-1)$, $m\in{\mathbb N}$, and ${\M}(L(m,m-1))$
is the disjoint union of ${\M}_1$ and ${\M}_2$ with
\[ {\M}_1=\{ a\in {\mathbb C}\co 0< \mbox{\rm Re}(a) < 1\} /
(a\sim 1-a),\]
\[ {\M}_2=\{ n\in {\mathbb N}\co n\equiv -1\;\mbox{\rm mod}\; m\} .\]
\end{thm}

\begin{rem}
{\rm Our notational convention is that ${\mathbb N}$ stands for the
natural numbers~$\geq 1$.}
\end{rem}

Of course, the component ${\M}_2$ of ${\M}(L(m,m-1))$ is
in any case simply a countable set, but the values $n\equiv -1$ mod~$m$
have concrete geometric meaning, namely, the contact circles corresponding
to $n\in {\M}_2$ are induced from the contact circle $(\omega_1,
\omega_2)$ on $S^3\subset {\mathbb C}^2$ defined by
\[ \omega_1+i\omega_2=nz_1dz_2-z_2dz_1+z_2^ndz_2.\]
The contact circles corresponding to $a\in{\M}_1$ are given by
\[ \omega_1+i\omega_2 =az_1dz_2-(1-a)z_2dz_1.\]
As pointed out in \cite{gego95} one cannot put a `good' metric topology
on the disjoint union ${\M}_1\sqcup {\M}_2$. An analogous
phenomenon was first observed by Kodaira-Spencer~\cite{kosp58}
in their study of moduli of Hopf surfaces.

The Teichm\"uller spaces for the spherical manifolds will be
determined in Section~\ref{section:s} (Theorems~\ref{thm:teich-s1}
and~\ref{thm:teich-s2}).

Henceforth, $M$ will denote a closed, orientable $3$-manifold of
the form $\Gamma\backslash{\G}$ with ${\G}$ equal to
$\widetilde{\SL}_2$ or $\widetilde{\E}_2$ and $\Gamma$ a discrete, cocompact
subgroup of~${\G}$.
Write $\pi =\pi_1(M)\cong\Gamma$
for the fundamental group of $M$. Consider the pair
$(\pi ,{\G})$, where ${\G}$ is the geometry on which $M$
is modelled. The first main result of the present paper translates
the geometric definition of Teichm\"uller and moduli spaces into
a more algebraic one. This will later allow us to arrive at explicit
descriptions of these deformation spaces. First we recall
the following concept, cf.\ \cite{harv77}, \cite{klr85}, \cite{vish93}.

\begin{defn}
The {\bf Weil space} ${\R}(\pi ,{\G})$ is
the space of faithful
representations $\rho$ of $\pi$ in ${\G}$ such that $\rho (\pi )$
is discrete and cocompact in~${\G}$.
\end{defn}

Elements $a\in {\G}$ act from the left on
${\R}(\pi ,{\G})$ by inner automorphisms
\[ \rho\longmapsto\rho^a,\;\; \rho^a(u)=a\rho(u)a^{-1}.\]
Write $\Inn ({\G})$ for the group of such inner automorphisms.
The automorphism group $\Aut (\pi )$ of $\pi$ acts from the right on
${\R}(\pi ,{\G})$; an element $\vartheta\in\Aut (\pi )$
acts by
\[ \rho\longmapsto\rho\circ \vartheta .\]
Write $\Out (\pi )$ for the group of outer automorphisms of~$\pi$.

In the case ${\G}=\widetilde{\E}_2$ we make a slight modification
to the Weil space. Recall that $\widetilde{\E}_2$ may be regarded
as ${\mathbb R}^3$ with multiplication
\begin{eqnarray*}
\lefteqn{
\left( \left( \begin{array}{c}x_0\\ y_0\end{array} \right) ,\theta_0\right)
\cdot \left( \left( \begin{array}{c} x\\ y\end{array} \right) ,
\theta \right) }\\
  & = &
\left( \left( \begin{array}{cc}\cos\theta_0 & -\sin\theta_0\\
 \sin\theta_0 & \cos\theta_0\end{array} \right)
\left( \begin{array}{c}x\\y\end{array} \right) +
\left( \begin{array}{c}x_0\\y_0\end{array} \right) ,
\theta_0 +\theta \right).
\end{eqnarray*}
We call the $(x,y)$-component of an element in $\widetilde{\E}_2$ its
{\bf translational part} and the $\theta$-component its
{\bf rotational part}. Define
\[ {\R}'(\pi ,\widetilde{\E}_2)={\R}(\pi ,\widetilde{\E}_2)
/{\mathbb R}^+,\]
where ${\mathbb R}^+$ acts on the Weil space by scaling of the translational
parts of the $\rho (u)$, $u\in\pi$.

\begin{defn}
\label{defn:algebraic}
The {\bf algebraic Teichm\"uller space} $\T^{\alg}(M)$ of $M$ is defined as
\[ \Inn (\widetilde{\SL}_2)\backslash {\R}(\pi, \widetilde{\SL}_2) \]
if $M$ is modelled on~$\widetilde{\SL}_2$, and as
\[ \Inn (\widetilde{\E}_2)\backslash {\R}'(\pi, \widetilde{\E}_2) \]
if $M$ is modelled on~$\widetilde{\E}_2$.
\end{defn}

\begin{thm}
\label{thm:deform-alg}
If $M$ is modelled on $\widetilde{\SL}_2$ or $\widetilde{\E}_2$, then
\[ \T (M)= \T^{\alg}(M) \]
and
\[ \M (M)=\T^{\alg}(M)/\Out (\pi ).\]
\end{thm}

Initially we shall treat the equal signs in this theorem merely
as set-theoretic bijections. In the course of our further explicit
description of these spaces we shall see that there are in fact natural
topologies on these deformation spaces and the various quotients of
the Weil space such that the equal signs may be read as homeomorphisms.

If $M$ is modelled on $\widetilde{\SL}_2$, it admits a Seifert fibration
$M\rightarrow O_M$, unique up to isotopy. Write $\T(O_M)$ for the Teichm\"uller
space of hyperbolic structures on the base orbifold~$O_M$, and write
$\T^{\alg}(O_M)$ for the algebraic Teichm\"uller space of $O_M$, defined
via representations of the orbifold fundamental group in~$\PSL_2{\mathbb R}$.
As will be explained in Section~\ref{section:sl},
a taut contact circle on $M$ induces in a natural way a hyperbolic structure
on~$O_M$.
We shall exhibit a commutative diagram
\[
\begin{CD}
\T (M) @>\Phi>> \T^{\alg}(M)\\
@VVV             @VVV\\
\T (O_M) @>\phi>> \T^{\alg}(O_M),
\end{CD}
\]
where $\Phi$ and $\phi$ are homeomorphisms, and the vertical maps
are coverings.

Passing to moduli spaces, our theory entails the construction of certain
non-trivial branched coverings of the moduli space $\M (O_M)$ of
hyperbolic structures on~$O_M$. Some comments about this issue are
made at the end of Section~\ref{section:sl}; we reserve details for
a forthcoming publication~\cite{gegob}.

The complex structure on the `classical' Teichm\"uller space
$\T(O_M)$ induces complex structures on the other spaces in this diagram
so that all maps become local biholomorphisms.
By exhibiting a universal family over~$\T (M)$
it will be shown in Section~\ref{section:complex} that
this complex structure on $\T (M)$ is the natural one to consider.
In the case of the geometries $S^3$ and $\widetilde{\E}_2$, Teichm\"uller
or moduli space is (or will be) described explicitly as a complex space.
In the $\widetilde{\E}_2$ case we also describe a universal family. In the
$S^3$ case it is clear from the explicit description of the taut contact
circles corresponding to $\M_1$ that the complex structure on $\M_1$
is naturally adapted to the classification of taut contact circles.
Further justification of this will also be given in
Section~\ref{section:complex}, where we describe a complex Godbillon-Vey
invariant, which turns out to define a holomorphic isomorphism
of $\M_1$ with a domain in~${\mathbb C}$.

Theorem~\ref{thm:deform-alg} shows that
taut contact circles on left-quotients of $\widetilde{\SL}_2$ or
$\widetilde{\E}_2$ have deformation spaces closely related to the
deformation spaces for geometric structures (in the sense of
Thurston) on Seifert fibred $3$-manifolds studied by
Kulkarni-Lee-Raymond~\cite{klr85}. For these geometric structures, the
translation from a geometric to an algebraic definition of the
corresponding deformation spaces is given by the developing map. This is not
made explicit in~\cite{klr85}, where the investigation starts directly
from the algebraic definition. In our situation, slightly more work is
necessary for this process of translation.
It needs to be pointed out
that our algebraic set-up differs from that of~\cite{klr85} by the
fact that we consider representations in ${\G}$ rather than in the
isometry group of~${\G}$. Moreover, our subsequent
elaboration of these deformation spaces in Sections \ref{section:sl}
and~\ref{section:e} is closer in spirit to the work of Ohshika~\cite{ohsh87}.
Compare also the related results of Ue~\cite{ue92} about deformations of
geometric structures on Seifert $4$-manifolds. Again, that paper deals
directly with the algebraic definitions of the relevant deformation spaces.

We should like to argue that the deformation spaces for taut contact
circles introduced in the present paper may prove valuable on two counts
-- apart from the fact mentioned above that they give rise to
interesting non-trivial branched coverings of~$\M (O_M)$.
First of all, they provide a geometric interpretation (in a slightly
altered setting) of the algebraic deformation spaces studied by
Kulkarni-Lee-Raymond, and thus carry structure not seen in~$\T^{\alg}(M)$. 
This additional structure might be of interest in getting a better
understanding of the related algebraic spaces, in particular the Weil space,
just as classical (geometric) Teichm\"uller theory, notably
the analytical methods developed there, yield information about
the Weil space of representations in $\PSL_2{\mathbb R}$.

Secondly, our deformation spaces may be seen as a bridge between the
deformation theory of 2-dimensional orbifolds and that of complex
surfaces. In the case of $\widetilde{\SL}_2$, it also makes sense
to regard our theory as some sort of desingularisation of the
2-dimensional theory.

Theorem~\ref{thm:deform-alg} will be proved in Section~\ref{section:Cartan}.
The explicit geometric descriptions of these deformation spaces -- the
analogues of Theorem~\ref{thm:deform-s} -- constitute the other principal
results of the present paper. For the geometry $\widetilde{\SL}_2$ this
description is given in Theorem~\ref{thm:deform-sl}; for $\widetilde{\E}_2$,
in Theorem~\ref{thm:deform-e}.

\begin{rem}
{\rm In our definition of a taut contact circle we have opted not to make
any assumptions on orientations. That is, we do not require that
$\omega_i\wedge d\omega_i$ be a positive volume form for a chosen
orientation on~$M$. As it turns out, most of the manifolds that admit
taut contact circles, e.g.\ left-quotients of $\widetilde{\SL}_2$ and
non-abelian left-quotients of~$\SU (2)$, do not admit any
orientation reversing diffeomorphisms, see~\cite{nera78}.
In these cases, it will be implicit in
our results that all taut contact circles on a given $M$ define one and the
same orientation.}
\end{rem}

By taking triples of $1$-forms in place of pairs, but otherwise making
requirements analogous to those in Definition~\ref{defn:circle},
we arrive at the notion of a {\bf taut contact sphere}. As shown
in~\cite{gego95}, such structures exist exactly on the left-quotients
of~$\SU (2)$. The corresponding moduli problem is solved in~\cite{gegoa},
by methods quite different from those in the present paper.
No higher-dimensional linear families of
contact forms can exist on $3$-manifolds, because our definition
implies in particular the pointwise linear independence of the $1$-forms
spanning the linear family. Some generalisations of these concepts to
higher-dimensional manifolds have been discussed in~\cite{gego96}
and~\cite{geth95}.

\section{Cartan structures and Weil spaces}
\label{section:Cartan}
In this section we prove Theorem~\ref{thm:deform-alg}. We begin with several
preliminary observations about Cartan structures.
The idea is to show that, in the cases of interest to us,
each conformal class of taut contact circles contains an essentially
unique distinguished Cartan structure (Proposition~\ref{prop:uniqueK}).
With the help of this
Cartan structure we can define an analogue of the developing map
for geometric structures (Lemma~\ref{lem:f}), which will be the key
to proving Theorem~\ref{thm:deform-alg}.

\begin{lem}
Let $(\omega_1,\omega_2)$ be a Cartan structure on some $3$-manifold. Then
there is a unique $1$-form $\omega_3$ such that
\begin{eqnarray*}
d\omega_1 & = & \omega_2\wedge\omega_3,\\
d\omega_2 & = & \omega_3\wedge\omega_1.
\end{eqnarray*}
\end{lem}

The proof of this simple lemma is left to the reader. Occasionally we write
a Cartan structure as $(\omega_1,\omega_2,\omega_3)$ if we wish to
fix a label for the third 1-form determined by $(\omega_1,\omega_2)$.

\begin{defn}
A Cartan structure $(\omega_1,\omega_2)$ is called a $K$-{\bf Cartan structure}
if the unique $\omega_3$ of the preceding lemma satisfies
\[ d\omega_3=K\omega_1\wedge\omega_2.\]
\end{defn}

Here $K$ may in principle be any function on~$M$ that is constant
along the common kernel of $\omega_1$ and $\omega_2$ (since
$dK\wedge\omega_1\wedge\omega_2=d^2\omega_3=0$), but the two
cases of interest
to us will be $K\equiv -1$ and $K\equiv 0$. Notice that the Lie groups
${\mathcal G}$ in Theorem~\ref{thm:classification} are exactly the
universal covers of the group of orientation preserving isometries
of the simply-connected $2$-dimensional space forms. For $K$ a constant,
write ${\mathcal G}_K$ for the corresponding Lie group.
We mention without proof the
following variation on Theorem~\ref{thm:classification},
which can be derived without much trouble from
the results of~\cite{gego95}.

\begin{prop}
Let $M$ be a closed $3$-manifold. Then $M$ admits a $K$-Cartan structure
with $K$ constant if and only if $M$ is diffeomorphic to a left-quotient
$\Gamma\backslash {\G}_K$.
\end{prop}

\begin{rem}
{\rm A $3$-manifold with a $(-1)$-Cartan structure is equivalent to what
Jacobowitz~\cite{jaco90} calls a {\bf projective structure}. These
structures are defined as a $3$-dimensional manifold together with a triple of
independent $1$-forms $\beta_1,\beta_2,\beta_3$ satisfying the equations
\begin{eqnarray*}
d\beta_1 & = & -\beta_2\wedge \beta_3,\\
d\beta_2 & = & -2\beta_1\wedge\beta_2\\
d\beta_3 & = & 2\beta_1\wedge\beta_3.
\end{eqnarray*}
The translation from a projective structure to a $(-1)$-Cartan structure is
given by $\omega_1=2\beta_1$, $\omega_2=\beta_2+\beta_3$, $\omega_3=
\beta_2-\beta_3$.}
\end{rem}

\begin{prop}
\label{prop:uniqueK}
{\rm (a)} Let $M$ be a (compact) left-quotient of $\widetilde{\SL}_2$.
In each conformal class of taut contact circles on $M$ there is one and only
one $(-1)$-Cartan structure.

{\rm (b)} Let $M$ be a (compact) left-quotient of $\widetilde{\E}_2$.
Every conformal
class of taut contact circles on $M$ contains a $0$-Cartan structure.
This $0$-Cartan structure is unique up to multiplication by a
positive constant.
\end{prop}

\noindent {\em Proof.}
(a) Existence of a $(-1)$-Cartan structure in every homothety
class was proved in~\cite[Thm.~7.4]{gego95}; since the rotate of a
$K$-Cartan structure is again a $K$-Cartan structure (with the same
$\omega_3$) we also have existence in every conformal class.
To prove uniqueness,
we argue as follows. Let $(\omega_1,\omega_2)$ be a $(-1)$-Cartan
structure on $M$. As shown in the cited theorem from~\cite{gego95}, the
common kernel $\ker\omega_1\cap\ker\omega_2$ determines a Seifert fibration
$M\rightarrow O_M$, with $O_M$ a hyperbolic orbifold.
Lift $(\omega_1,\omega_2)$ to a finite covering
$M'\rightarrow M$ corresponding to a covering $\Sigma_g\rightarrow O_M$
of $O_M$ by an honest surface.
The $(-1)$-Cartan equations imply that the symmetric
bilinear form $\omega_1\otimes\omega_1+\omega_2\otimes\omega_2$ is 
invariant under the flow of $X\in\ker\omega_1\cap\ker\omega_2$ with
$\omega_3(X)=1$, and that it induces a metric of constant curvature~$-1$
on~$\Sigma_g$.

Now assume that we have a further $(-1)$-Cartan structure $(v\omega_1,
v\omega_2)$ on~$M$. The structure equations imply $dv\wedge\omega_1\wedge
\omega_2=0$, so $v$ is constant along the fibres of the $S^1$-fibration,
and $v^2(\omega_1\otimes\omega_1+\omega_2\otimes\omega_2)$ also
induces a metric of constant curvature~$-1$ on~$\Sigma_g$.

But there is a unique hyperbolic metric (of constant curvature~$-1$)
in a given conformal class of metrics on~$\Sigma_g$, since
conformal automorphisms of the hyperbolic plane ${\mathbb H}^2$
are actually isometries. This forces $v\equiv 1$.

\vspace{2mm}

(b) In the euclidean case, existence of a $0$-Cartan structure in every
conformal class was also proved in~\cite{gego95}. However, the common
kernel of $\omega_1$ and $\omega_2$ will not, in general, define a Seifert
fibration (cf.~\cite[Thm.~7.4]{gego95}),
so we need a substitute for the argument in~(a).

Assume that $(\omega_1,\omega_2,\omega_3)$ and $(v\omega_1,v\omega_2,
\omega_3')$ are $0$-Cartan structures on a given left-quotient
of~$\widetilde{\E}_2$. As in (a) the structure equations imply
$dv\wedge\omega_1\wedge\omega_2=0$. Since
\[ \omega_1\wedge\omega_2\wedge\omega_3=\omega_1\wedge d\omega_1\neq 0,\]
the $\omega_i$ are pointwise linearly independent and we can write
\[ dv=v_1\omega_1+v_2\omega_2\]
with uniquely defined functions $v_i\co M\rightarrow {\mathbb R}$.
Then the structure equations for the two $0$-Cartan structures yield
\begin{eqnarray*}
v_2\omega_2\wedge\omega_1+v\omega_2\wedge\omega_3 & = & v\omega_2
                              \wedge\omega_3',\\
v_1\omega_1\wedge\omega_2+v\omega_3\wedge\omega_1 & = &
                   v\omega_3'\wedge\omega_1.
\end{eqnarray*}
This implies
\[ \omega_3'=\frac{v_2}{v}\omega_1-\frac{v_1}{v}\omega_2+\omega_3\]
and
\[ d(\frac{v_2}{v}\omega_1-\frac{v_1}{v}\omega_2)=d(\frac{v_2}{v}
\omega_1-\frac{v_1}{v}\omega_2+\omega_3)
     =d\omega_3'=0.\]
Let $\ast$ be the Hodge star with respect to the Riemannian metric
\[ \omega_1\otimes\omega_1+\omega_2\otimes\omega_2+\omega_3\otimes\omega_3\]
on $M$ (and orientation given by $\omega_1\wedge\omega_2\wedge\omega_3$).
Then
\[ \ast d\log v=\frac{\ast dv}{v}=
\frac{v_1}{v}\omega_2\wedge\omega_3+\frac{v_2}{v}\omega_3\wedge\omega_1=
\frac{1}{v} (v_1\omega_2-v_2\omega_1)\wedge\omega_3,\]
and hence $d\!\ast\! d\log v=0$.
On functions on 3-manifolds, the Laplace operator
$\Delta$ takes the form $\Delta =-\ast\! d\!\ast\! d$, so we have in fact
$\Delta \log v=0$, i.e.\ $\log v$ is harmonic. Since $M$ is closed, this
forces $\log v$ and hence $v$ to be constant. \hfill $\Box$

\begin{rem}
{\rm There is no analogue of this result for left-quotients of~$\SU (2)$,
see~\cite[Section~6]{gego95}. It is for this reason that in the
spherical case there is no direct translation into algebraic
Teichm\"uller and moduli spaces.}
\end{rem}

Fix a standard left-invariant $K$-Cartan structure $(\alpha_1,\alpha_2,
\alpha_3)$ on ${\G}$, with $({\G},K)$ equal to
$(\widetilde{\SL}_2,-1)$ or $(\widetilde{\E}_2,0)$,
cf.~\cite[Section~2]{gego95}. In particular, on $\widetilde{\E}_2$ we may
choose
\begin{eqnarray*}
\alpha_1 & = & \cos\theta\, dx+\sin\theta\, dy,\\
\alpha_2 & = & -\sin\theta\, dx+\cos\theta\, dy,\\
\alpha_3 & = & -d\theta
\end{eqnarray*}
For each discrete, cocompact subgroup $\Gamma\subset\G$, we continue
to write $\alpha_i$ for the $1$-forms induced on~$\Gamma\setminus
\G$. If $M$ is a compact manifold diffeomorphic to
a left-quotient of~${\G}$, and $(\omega_1,\omega_2,
\omega_3)$ a $K$-Cartan structure on~$M$, we keep the same notation for the
pull-back of this $K$-Cartan structure to the universal cover~$\widetilde{M}$.

\begin{lem}
\label{lem:f}
There is a diffeomorphism $f\co\widetilde{M}\rightarrow{\cal G}$
such that $f^*(\alpha_i)=\omega_i$, $i=1,2,3$, and such that $f$
descends to
a diffeomorphism $\overline{f}\co M\rightarrow\Gamma\backslash
{\G}$, for a suitable discrete, cocompact subgroup $\Gamma$
of~${\G}$.
\end{lem}

\begin{rem}
{\rm This diffeomorphism $f$ is the analogue of the developing map
for $(G,X)$-structures in the sense of Thurston~\cite{thur97}.}
\end{rem}

\noindent {\em Proof.}
According to~\cite[Thm.~1.6]{gego95}, given any taut contact circle
$(\omega_1,\omega_2)$ on~$M$, there is a $\Gamma\subset {\G}$
and a diffeomorphism $\overline{f}\co M\rightarrow \Gamma\backslash
{\G}$ which pulls back the standard Cartan structure $(\alpha_1,
\alpha_2)$ on $\Gamma\backslash {\G}$ to the homothety class
of $(\omega_1,\omega_2)$ on~$M$. Rotations of $(\alpha_1,\alpha_2)$
on~${\G}$ (without changing~$\alpha_3$) can be effected by the flow
of the Lie algebra element $e_3$ dual to~$\alpha_3$ (with respect to the
basis $\alpha_1,\alpha_2,\alpha_3$). This is the same as right
multiplication by the one-parameter subgroup corresponding to~$e_3$,
and so it descends to any left-quotient $\Gamma\backslash {\G}$.
This implies that we can actually find a diffeomorphism $\overline{f}\co
M\rightarrow \Gamma\backslash {\G}$ that pulls back $(\alpha_1,
\alpha_2)$ to the {\em conformal class} of~$(\omega_1,\omega_2)$.

Now both $(\omega_1,\omega_2)$ and $\overline{f}^*(\alpha_1,\alpha_2)$
are $K$-Cartan structures (the former by assumption, the latter by
construction). For $({\G},K)=(\widetilde{\SL}_2,-1)$ this
forces $(\omega_1,\omega_2)=\overline{f}^*(\alpha_1,\alpha_2)$ by
Proposition~\ref{prop:uniqueK}; for $({\G},K)=
(\widetilde{\E}_2,0)$ we only get $(\omega_1,\omega_2)= c\overline{f}^*
(\alpha_1,\alpha_2)$ for some $c\in{\mathbb R}^+$.

Consider the Lie group isomorphism
\[ \begin{array}{rccc}
s_c\co & \widetilde{\E}_2 & \longrightarrow & \widetilde{\E}_2\\
       & \left( \left(\begin{array}{c}x\\y\end{array}\right) ,\theta
         \right) & \longmapsto &
         \left( c \left(\begin{array}{c}x\\y\end{array}\right) ,\theta
         \right) .
\end{array} \]
Then $s_c^*(\alpha_1,\alpha_2,\alpha_3)=(c\alpha_1,c\alpha_2,\alpha_3)$.
Hence we obtain the
desired diffeomorphism if we replace the lift $f$ of
$\overline{f}\co M\rightarrow\Gamma\backslash {\G}$ by $s_c\circ f$,
and the subgroup $\Gamma$ by $s_c(\Gamma )$. \hfill $\Box$

\vspace{2mm}

This diffeomorphism $f$ will determine a representation of $\pi =\pi_1(M)$
in~${\G}$. The following lemma describes the indeterminacy in this
construction.

\begin{lem}
\label{lem:choice}
Let $f$ and $f'$ be two diffeomorphisms $\widetilde{M}\rightarrow {\G}$
as in the preceding lemma, with corresponding subgroups $\Gamma ,
\Gamma '\subset\G$. Then there is a $w\in {\G}$ such that
$f'=L_w\circ f$ (and $\Gamma '=w\Gamma w^{-1}$), with $L_w$ denoting
left multiplication.
\end{lem}

\noindent {\em Proof.}
Fix a base point $p_0\in\widetilde{M}$. There is a $w\in{\G}$ such
that $L_w\circ f(p_0)=f'(p_0)$.
The triple of $1$-forms $\omega_1-\alpha_1$, $\omega_2-\alpha_2$,
$\omega_3-\alpha_3$ on $\widetilde{M}\times {\G}$ (where we
forgo writing the pull-back maps of the projections onto the two
factors) generate a differential ideal ${\mathcal I}$, that is, $d{\mathcal I}
\subset {\mathcal I}$. For instance,
\[ d(\omega_1-\alpha_1) = \omega_2\wedge (\omega_3-\alpha_3)+
                         (\omega_2-\alpha_2)\wedge\alpha_3.\]
So ${\mathcal I}$ defines a $3$-dimensional foliation on $\widetilde{M}
\times {\G}$. The graphs of $f'$ and $L_w\circ f$ are integral
submanifolds of this foliation, and since both pass through the
point $(p_0,f'(p_0))=(p_0,L_w\circ f(p_0))$, they must be
identical.\hfill $\Box$

\vspace{2mm}

\noindent {\em Proof of Theorem}~\ref{thm:deform-alg}.
The key to the proof of this theorem is that we have the following
ingredients:
\begin{itemize}
\item[(1)] a developing map (constructed in the preceding lemmas),
which allows to translate from taut contact circles to suitable
representations of the fundamental group;
\item[(2)] an analogue of the Nielsen theorem, which guarantees that
all automorphisms of the fundamental group are induced by diffeomorphisms,
and thus ensures the surjectivity of this process of translation;
\item[(3)] an analogue of the Baer theorem, which allows to pass from
homotopy to isotopy, and thus ensures injectivity.
\end{itemize}
Any geometric structure (on manifolds homotopy equivalent to
Eilenberg-MacLane spaces) that provides these ingredients will
have an algebraic deformation theory based on representations
of the fundamental group analogous to the one we are about to
construct.

\vspace{2mm}

Identify $\pi$ with the deck transformation group of the universal
covering $\widetilde{M}\rightarrow M$. Thus we regard the elements of
$\pi$ as acting on $\widetilde{M}$ from the left.

Given a taut contact circle on~$M$, choose a $K$-Cartan structure in its
homothety class. Choose a diffeomorphism $f\co\widetilde{M}\rightarrow
{\G}$, lifted from a diffeomorphism $\overline{f}\co M\rightarrow
\Gamma\backslash{\G}$, with $f^*(\alpha_i)=\omega_i$, $i=1,2,3$.

Given a deck transformation $u\in\pi$ of $\widetilde{M}\rightarrow M$, we get
a corresponding deck transformation
\[ f\circ u\circ f^{-1}\]
of ${\G}\rightarrow \Gamma\backslash{\G}$. This must of
course be an element of $\Gamma\subset{\G}$. Indeed, $f\circ u\circ
f^{-1}$ acts trivially on a basis $(\alpha_1,\alpha_2,\alpha_3)$ of
left-invariant $1$-forms, and hence on every left-invariant $1$-form;
therefore $f\circ u\circ f^{-1}$ is left multiplication by some
element $\rho (u)\in {\G}$.

So a choice of $K$-Cartan structure and a choice of $f$ determine a
discrete, faithful representation $\rho\in{\R}(\pi,{\G})$
via the equation
\[ f\circ u=L_{\rho (u)}\circ f.\]
The following points (i)--(iii) discuss the dependence of
$\rho$ on the various choices in this
construction.

\vspace{2mm}

(i) According to Lemma~\ref{lem:choice}, a different choice $f'$ instead of
$f$ (with fixed $K$-Cartan structure) must be of the form $f'=L_w\circ f$
for some $w\in{\G}$. The corresponding representation $\rho '$
satisfies
\[ f'\circ u=L_{\rho '(u)}\circ f'\]
for all $u\in\pi$, that is,
\[ L_w\circ f\circ u=L_{\rho '(u)}\circ L_w\circ f.\]
Hence
\[ L_{\rho (u)}\circ f=f\circ u=L_w^{-1}\circ L_{\rho '(u)}\circ
L_w\circ f=L_{w^{-1}\rho '(u)w}\circ f,\]
and we conclude $w^{-1}\rho '(u)w=\rho (u)$. So different choices of $f$ lead
to representations that differ by an element of $\Inn ({\G})$.

\vspace{2mm}

(ii) We can change a $K$-Cartan structure $(\omega_1,\omega_2,\omega_3)$
on $M$ within its homothety class by rotating $(\omega_1,\omega_2)$ by an
angle $\theta$ and keeping $\omega_3$ fixed. Write $R_{\theta}$ for the
time~$\theta$ flow of $e_3$ on~${\G}$ (in the notation of the proof
of Lemma~\ref{lem:f}). The map $R_{\theta}$ is given by right multiplication,
and it has the effect of keeping $\alpha_3$ fixed and pulling back
$(\alpha_1,\alpha_2)$ to a $\theta$-rotate of itself.

So replacing $(\omega_1,\omega_2)$ by a $\theta$-rotate amounts to
replacing $f$ by $f'=R_{\theta}\circ f$. From the defining equation for
the corresponding representation~$\rho '$,
\[ R_{\theta}\circ f\circ u=L_{\rho '(u)}\circ R_{\theta}\circ f,\]
we get
\[ L_{\rho (u)}=f\circ u\circ f^{-1}=R_{\theta}^{-1}\circ L_{\rho '(u)}
\circ R_{\theta}.\]
But left and right multiplication commute, so $\rho '(u)=\rho (u)$.

\vspace{2mm}

(iii) By Proposition~\ref{prop:uniqueK}, the rotation discussed in (ii)
is the only ambiguity in the choice of a $(-1)$-Cartan structure within
a homothety class of taut contact circles on a left-quotient
of~$\widetilde{\SL}_2$. For left-quotients of $\widetilde{\E}_2$ we may
also scale the $0$-Cartan structure by a positive constant. If
$f\co\widetilde{M}\rightarrow\widetilde{\E}_2$ pulls back $(\alpha_1,
\alpha_2,\alpha_3)$ to $(\omega_1,\omega_2,\omega_3)$, then
$f'=s_c\circ f$, $c\in{\mathbb R}^+$, pulls it back to
$(c\omega_1,c\omega_2,\omega_3)$. So the representation $\rho '$ corresponding
to $(c\omega_1,c\omega_2)$ is defined by
\[ s_c\circ f\circ u=L_{\rho '(u)}\circ s_c\circ f,\]
hence
\[ L_{\rho '(u)}=s_c\circ L_{\rho (u)}\circ s_c^{-1}.\]
If
\[ \rho (u)=\left(\left(\begin{array}{c}x_u\\y_u\end{array}\right) ,\theta_u
\right) ,\]
then a straightforward calculation gives
\[ \rho '(u)=\left( c \left(\begin{array}{c}x_u\\y_u\end{array}\right) ,
\theta_u \right) ,\]
so $\rho$ and $\rho'$ define the same element of $\R '(\pi ,
\widetilde{\E}_2)$.

\vspace{2mm}

From (i)--(iii) we conclude that there is a well-defined map
\[ \widetilde{\Phi}\co{\mathcal C}(M)\longrightarrow
{\T}^{\mbox{\scriptsize\rm alg}} (M), \]
where ${\T}^{\mbox{\scriptsize\rm alg}}$ is the algebraic
Teichm\"uller space introduced in Definition~\ref{defn:algebraic}.

\vspace{2mm}

(iv) The next step is to see that $\widetilde{\Phi}$ induces a map
\[ \Phi\co {\T}(M)\longrightarrow
{\T}^{\mbox{\scriptsize\rm alg}}(M).\]
To this end, suppose that we take the pull-back of a given taut contact
circle $(\omega_1,\omega_2)$ on $M$ under $\overline{\varphi}\in
\Diff_0(M)$. The corresponding $K$-Cartan structures in the respective
homothety classes may also be assumed to be related by the pull-back
map~$\overline{\varphi}^*$. The diffeomorphism $\overline{\varphi}$ lifts
to a diffeomorphism $\varphi$ of~$\widetilde{M}$. Then
$f\circ\varphi\co\widetilde{M}\rightarrow {\G}$ is a lift
of $\overline{f}\circ\overline{\varphi}\co M\rightarrow \Gamma\backslash
{\G}$ and $(f\circ\varphi )^*\alpha_i=\varphi^*\omega_i$,
$i=1,2,3$. So the representation $\rho '$ corresponding to
$\overline{\varphi}^*
(\omega_1,\omega_2)$ is defined by
\[ f\circ\varphi\circ u=L_{\rho '(u)}\circ f\circ\varphi .\]
But $\varphi\circ u\circ \varphi^{-1}=u$, because $\oph$ is homotopic
to the identity, so $\rho '(u)=\rho (u)$.

\vspace{2mm}

(v) To complete the proof of the statement in Theorem~\ref{thm:deform-alg}
about Teich\-m\"uller spaces, we need to show that $\Phi$ is a bijection.

\vspace{2mm}

$\Phi$ is surjective: To prove this we first observe that the diffeomorphism
type of $\rho (\pi )\backslash {\G}$ is completely determined
by the isomorphism type of~$\pi$: For ${\G}=\widetilde{\SL}_2$ this
is a consequence of the fact that
$\rho (\pi )\backslash\widetilde{\SL}_2$
is a large Seifert manifold in the sense of Orlik~\cite[p.~92]{orli72}
with fundamental group~$\pi$,
see~\cite[p.~97]{orli72}. For ${\G}=\widetilde{\E}_2$ this
follows from the complete classification of left-quotients of
$\widetilde{\E}_2$; up to diffeomorphism these are exactly the five torus
bundles over $S^1$ with periodic monodromy, and it is easy to
compute explicitly that they are distinguished by their fundamental
groups, see Section~\ref{section:e}.

Therefore any $\rho\in\R (\pi ,\G )$ gives rise to a diffeomorphic
copy $\rho (\pi )\backslash\G$ of~$M$. Moreover, the diffeomorphism
$M\rightarrow \rho (\pi )\backslash\G$ can be chosen such that its induced
map on fundamental groups gives the representation~$\rho$, because
for the manifolds under consideration any automorphism of the fundamental
group can be induced by a diffeomorphism of the manifold.

For the
left-quotients of $\widetilde{\E}_2$ this last fact was proved by
Neuwirth~\cite{neuw63}.
For the left-quotients of $\widetilde{\SL}_2$
we argue as follows. First one observes that in this case $\pi$ contains
a unique maximal normal subgroup isomorphic to~${\mathbb Z}$,
generated by the class of a regular fibre of the Seifert fibration
$M\rightarrow O_M$, cf.~\cite{scot83} and Section~4 below. Thus an
automorphism of $\pi$ induces an automorphism of the orbifold
fundamental group $\pi^{\orb}=\pi_1^{\orb}(O_M)$, i.e.\ the deck
transformation group of the universal covering of~$O_M$ (which
in this case is the hyperbolic plane~${\mathbb H}^2$). By the
generalised Nielsen theorem, cf.~\cite[Thm.~8.1]{zies73},
this automorphism of $\pi^{\orb}$ is induced by conjugation
with a diffeomorphism $\psi$ of~${\mathbb H}^2$. This $\psi$
has to send any elliptic element of $\pi^{\orb}$ to another
elliptic element of the same order. In other words, $\psi$
descends to an `orbifold diffeomorphism' $\overline{\psi}$
of $O_M$ which preserves the set of cone points and may
permute only cone points of the same order amongst each other.
Moreover, the multiple fibres over cone points that are being
permuted must have the same Seifert invariants~$(\alpha ,\beta )$,
since we started with an automorphism
of~$\pi$. Therefore $\overline{\psi}$ lifts to a diffeomorphism
of~$M$, inducing the given automorphism of~$\pi$.

\begin{rem}
{\rm The preceding argument about lifting a diffeomorphism
$\overline{\psi}$ of $O_M$ to a diffeomorphism $\psi$
of~$M$ -- provided $\overline{\psi}$ corresponds to an
automorphism of $\pi^{\orb}$ induced from one of~$\pi$ --
is valid for all Seifert manifolds
fibred over a hyperbolic orbifold. For left-quotients~$M$
of $\widetilde{\SL}_2$ we shall in fact see in Section~\ref{section:sl}
that -- on a fixed~$M$ --  multiple fibres of the same order $\alpha$
always have the same Seifert invariant~$(\alpha ,\beta )$. This is not
true, in general, for quotients of $\widetilde{\SL}_2$ under discrete,
cocompact subgroups of the full isometry group of~$\widetilde{\SL}_2$,
or Seifert manifolds modelled on the geometry ${\mathbb H}^2\times
{\mathbb R}$. As a consequence, for left-quotients $M$
of $\widetilde{\SL}_2$ {\em any} diffeomorphism of $O_M$ lifts
to one of~$M$. This causes
some essential simplification as compared with the results
of~\cite{klr85}; see in particular the proof of
Lemma~\ref{lem:aut} below.}
\end{rem}

$\Phi$ is injective: Suppose we have two elements in ${\T}(M)$
with the same image in ${\T}^{\mbox{\scriptsize\rm alg}}(M)$
under~$\Phi$. By (i) and (iii) above we may choose $K$-Cartan
structures $(\omega_1,\omega_2)$ and $(\omega_1',\omega_2')$
representing these two elements and corresponding diffeomorphisms
$f,f'\co \widetilde{M}\rightarrow{\G}$ in such a way that
we get the same image even in ${\R}(\pi ,{\G})$, that is,
we may assume that for all $u\in\pi$ we have
\begin{eqnarray*}
f\circ u & = & L_{\rho (u)}\circ f,\\
f'\circ u & = & L_{\rho (u)}\circ f'.
\end{eqnarray*}
Then
\[ f^{-1}\circ f'\circ u=f^{-1}\circ L_{\rho (u)}\circ f'=
   u\circ f^{-1}\circ f',\]
which implies that $f^{-1}\circ f'$ descends to a diffeomorphism $\oph$
of $M$ which induces the identity on $\pi=\pi_1(M)$ and pulls back
$(\omega_1,\omega_2)$ to $(\omega_1',\omega_2')$. Since $M$ is an
Eilenberg-MacLane space $K(\pi ,1)$, this implies that $\oph$
is homotopic to the identity. By results of Waldhausen~\cite{wald68},
Scott~\cite{scot85}, and Boileau-Otal~\cite{boot91},
for the spaces in question this implies that
$\oph$ is isotopic to the identity. Thus $(\omega_1,\omega_2)$ and
$(\omega_1',\omega_2')$ define the same element in~${\T}(M)$.

\vspace{2mm}

(vi) Next we want to prove the statement of Theorem~\ref{thm:deform-alg}
concerning the moduli space~${\M}(M)$. Define
\[ \M^{\alg}(M)=\T^{\alg}(M)/\Out (\pi ).\]
So we have a map
\[ \T (M)\longrightarrow \M^{\alg}(M),\]
given by composing $\Phi$ with the quotient map $\T^{\alg}(M)\rightarrow
\M^{\alg}(M)$, and analogous to point (iv) above we need to check that this
descends to a map
\[ \Phi '\co\M (M)\longrightarrow\M^{\alg}(M).\]
Consider a taut contact circle $(\omega_1,\omega_2)$ on $M$ and its
pull-back under $\overline{\varphi}\in\Diff (M)$. As in (iv) we find for
the representations $\rho$ and $\rho '$ corresponding to $(\omega_1,\omega_2)$
and $\overline{\varphi}^*(\omega_1,\omega_2)$, respectively, that
\[ L_{\rho '(u)}\circ f=f\circ\varphi\circ u\circ\varphi^{-1}=
   L_{\rho (\varphi\circ u\circ\varphi^{-1})}\circ f.\]
Hence $\rho '(u)=\rho\circ \vartheta (u)$, where
$\vartheta\in\Aut (\pi)$ is defined by
$u\mapsto \varphi\circ u\circ\varphi^{-1}$, and hence $[\rho ']=[\rho ]$
in~$\M^{\alg}(M)$.

In spite of its appearance, the automorphism
$\vartheta$ is of course, in general,
not an inner automorphism of~$\pi$. Moreover, inner automorphisms of~$\pi$
do not move points in $\T^{\alg}(M)$, so it is appropriate to define
$\M^{\alg}(M)$ as the quotient of $\T^{\alg}(M)$ under $\Out (M)$.

\vspace{2mm}

(vii) Finally it remains to be checked that $\Phi '$ is bijective.
Surjectivity is proved as in (v) (without having to worry about
choosing a particular diffeomorphism $M\rightarrow\rho (\pi)\backslash\G$).
To prove injectivity, suppose that we have two elements in $\M (M)$ with
the same image in $\M^{\alg}(M)$ under~$\Phi '$. By (i) and (iii) we may
assume that these elements in $\M (M)$ are represented by $K$-Cartan
structures $(\omega_1,\omega_2)$ and $(\omega_1',\omega_2')$ which
already have the same image in
\[ {\mathcal S}(\pi )=\R (\pi, \G )/\Aut (\pi ),\]
which is usually called the {\bf Chabauty space} of~$\pi$ (or space
of discrete subgroups). So we have representations $\rho ,\rho '\in
\R (\pi ,\G )$ with $\rho '=\rho\circ \vartheta$ for some
$\vartheta\in\Aut (\pi )$, and
diffeomorphisms $f,f'\co\widetilde{M}\rightarrow\G$ that pull back
$(\alpha_1,\alpha_2)$ to $(\omega_1,\omega_2)$ and $(\omega_1',
\omega_2')$, respectively, such that
\begin{eqnarray*}
f\circ u & = & L_{\rho (u)}\circ f,\\
f'\circ u & = & L_{\rho\circ \vartheta (u)}\circ f'.
\end{eqnarray*}\
It is important to observe that both $f$ and $f'$ are lifted from
diffeomorphisms $M\rightarrow \Gamma\backslash\G$ with $\Gamma =\rho (\pi )
=\rho '(\pi )$. So $f^{-1}\circ f'$, as in~(v), descends to a diffeomorphism
$\oph$ of $M$ that pulls back $(\omega_1,\omega_2)$ to
$(\omega_1',\omega_2')$, so these $K$-Cartan structures define the same
element in~$\M (M)$. The automorphism $\vartheta$ of $\pi$ is determined by
the action of $\oph$ on~$\pi$, since
\[ f^{-1}\circ f'\circ u=\vartheta (u)\circ f^{-1}\circ f'.\]
This concludes the proof of Theorem~\ref{thm:deform-alg}.\hfill $\Box$
\section{$\widetilde{\SL}_2$-geometry}
\label{section:sl}
Let $M$ be a left-quotient of $\widetilde{\SL}_2$. Then $M$ is Seifert fibred,
and -- in contrast with the other geometries we are considering --
the Seifert fibration is unique up to isotopy, cf.\ \cite{scot83},
\cite[Cor.~2.3]{ohsh87}. In particular, the (closed, orientable) base
orbifold $O_M$ is uniquely determined by~$M$ (up to diffeomorphism).
Hence, if $(\omega_1,\omega_2)$ is a taut contact circle, its common
kernel will induce this unique Seifert fibration (cf.\ the proof of
Proposition~\ref{prop:uniqueK}). Moreover, if $(\omega_1,
\omega_2)$ is the $(-1)$-Cartan representative in that conformal
class, $\omega_1\otimes\omega_1 +\omega_2\otimes\omega_2$
will define a hyperbolic metric on~$O_M$,
and the coframe $(\omega_1,\omega_2)$ will induce an orientation
on~$O_M$. The pair $(\omega_2,\omega_1)$ would induce the same metric,
but the opposite orientation. So it would seem reasonable to fix
an orientation of $O_M$ once and for all, and only to consider taut
contact circles on $M$ compatible with this orientation. Yet, for
the translation to an algebraic deformation theory it is more
convenient not to fix such an orientation.

Thus, write $\C (O_M)$ for the space of hyperbolic metrics (of constant
curvature~$-1$) on $O_M$, together with a choice of orientation.
Equivalently, $\C (O_M)$ may be regarded as the space of
complex structures on~$O_M$. Teichm\"uller
and moduli space of such metrics are defined as
\begin{eqnarray*}
\T (O_M) & = & \C (O_M)/\Diff_0(O_M),\\
\M (O_M) & = & \C (O_M)/\Diff (O_M),
\end{eqnarray*}
These spaces come equipped with a natural topology, cf.~\cite{thur97}.
Notice that our definition
entails that $\T (O_M)$ has two components (corresponding
to the two choices of orientation), each homeomorphic
to ${\mathbb R}^{6g-6+2n}$, where $g$ is the genus
and $n$ the number of cone points of~$O_M$.

Again there are algebraic definitions of these deformation spaces. Write
\[ {\mathbb H}^2=\{\tau\in {\mathbb C}\co\mbox{\rm Im}\,\tau >0\} \]
for the upper half-plane. Its group of orientation preserving isometries
is $\PSL_2{\mathbb R}$, acting by fractional linear transformations.
As earlier
let $\pi^{\orb}$ denote the orbifold fundamental group $\pi_1^{\orb}(O_M)$,
i.e.\ the deck transformation group of the universal covering of~$O_M$
(an explicit presentation of this group will be given
shortly). Consider the Weil space
$\R (\pi^{\orb},\PSL_2{\mathbb R})$ of faithful representations
$\pi^{\orb}\rightarrow \PSL_2{\mathbb R}$ with discrete, cocompact image,
and set
\begin{eqnarray*}
\T^{\alg}(O_M) & = & \Inn (\PSL_2{\mathbb R})\backslash \R(\pi^{\orb},
\PSL_2{\mathbb R}),\\
\M^{\alg}(O_M) & = & \T^{\alg}(O_M)/\Out (\pi^{\orb}).
\end{eqnarray*}
Notice that $\T^{\alg}(O_M)$ has two components, just as $\T (O_M)$.
To get only one component, one would have to restrict to representations
that are related to one another by conjugation by an orientation preserving
homeomorphism of~${\mathbb H}^2$. This, however, would entail a not
entirely obvious restriction of the allowable automorphisms
of~$\pi^{\orb}$.

\begin{rem}
{\rm Let $\Pi$ be an abstract group topologised discretely, and $G$ a connected
Lie group (or at least a group such that the index of the identity
component $G_0$ in $G$ is finite). Write $\Hom (\Pi ,G)$ for the
space of representations $\Pi\rightarrow G$, equipped with the topology of
pointwise convergence. A fundamental result of Weil,
cf.~\cite{harv77}, \cite[p.~191]{vish93}, says that the Weil space
$\R (\Pi ,G)$ is open in
$\Hom (\Pi ,G)$. In the sequel, the algebraic Teichm\"uller and moduli spaces
will be equipped with the quotient topology coming from the topology
of $\R\subset\Hom$.}
\end{rem}

In classical Teichm\"uller theory it is shown that the maps
\[ \phi\co\T (O_M)\longrightarrow\T^{\alg}(O_M)\]
and
\[ \phi '\co\M (O_M)\longrightarrow\M^{\alg}(O_M)\]
given by the developing map are homeomorphisms, cf.~\cite[p.~301]{harv77},
 \cite[p.~194]{vish93}. In Section~\ref{section:Cartan}
we employed arguments analogous to these classical ones. The key
points in the classical theory is an analogue of Nielsen's
(and Baer's) theorem for orbifolds, establishing an equivalence
between isotopy classes of self-diffeomorphisms and outer
automorphisms of~$\pi^{\orb}$
(cf.~step~(v) in the proof of Theorem~\ref{thm:deform-alg}),
see~\cite{maha75}.

We recall a few further facts about Seifert manifolds from \cite{orli72},
\cite{scot83}. Assume that $M$ has normalised Seifert invariants
\[ \{ g,b,(\alpha_1,\beta_1),...,(\alpha_n,\beta_n)\} .\]
Since we are dealing with a left-quotient of $\widetilde{\SL}_2$, we know
that the orbifold Euler characteristic of the base
\[ \chi^{\orb}(O_M)=2-2g-n+\sum_{j=1}^n\frac{1}{\alpha_j}\]
and Euler number of the Seifert fibration
\[ e=-\Bigl( b+\sum_{j=1}^n\frac{\beta_j}{\alpha_j}\Bigr)\]
satisfy
\[ \chi^{\orb}(O_M)<0\;\;\mbox{\rm and}\;\; e\neq 0.\]
Furthermore, by a theorem of Raymond and Vasquez~\cite{rava81},
cf.~\cite{gego95a}, there are integers $r,k_1,...,k_n$ such that
\begin{eqnarray*}
rb & = & 2g-2-\sum_{j=1}^nk_j,\\
r\beta_j & = & \alpha_j-1+k_j\alpha_j,\;\; j=1,...,n.
\end{eqnarray*}
Observe that the integer $r$ is completely determined by the normalised Seifert
invariants, thanks to the formula $r=\chi^{\orb}/e$,
cf.\ Remark~\ref{rem:r} below.

\begin{rem}
\label{rem:mult}
{\rm The $\beta_j$ satisfy, by the definition of normalised Seifert
invariants, $1\leq \beta_j<\alpha_j$, and so
the equations above show that the multiple fibres with the same $\alpha_j$
also have the same~$\beta_j$. This is false, in general, for Seifert
manifolds which are not diffeomorphic to a left-quotient
of~$\widetilde{\SL}_2$.}
\end{rem}

In the sequel we fix the following presentations of $\pi^{\orb}
=\pi_1^{\orb}(O_M)$ and $\pi =\pi_1(M)$:
\begin{eqnarray*}
\pi^{\orb} & = & \Bigl\{ \ou_1,\ov_1,...,\ou_g,\ov_g,\oq_1,...,
\oq_n\co \prod_i[\ou_i,\ov_i]\prod_j\oq_j=1,\,\oq_j^{\alpha_j}=1\Bigr\} \\
\pi & = & \Bigl\{ u_1,v_1,...,u_g,v_g,q_1,...,q_n,h\co
\prod_i[u_i,v_i]\prod_jq_j=h^b,\\
  &   & \mbox{}\;\;\;\;\;\;\;\;q_j^{\alpha_j}h^{\beta_j}=1,\, h\;
\mbox{\rm central}\Bigr\}.
\end{eqnarray*}

Most of the work towards achieving explicit descriptions of the
deformation spaces in the $\widetilde{\SL}_2$-case will be contained
in the proof of the following theorem about Weil spaces.
Compare~\cite[Thm~2.5]{klr85} for the corresponding statement when
representations in $\Isom_0(\widetilde{\SL})$ are being considered.

\begin{thm}
\label{thm:weil-sl}
Let $M$ be a left-quotient of $\widetilde{\SL}_2$ with fundamental
group~$\pi$.
The Weil space $\R (\pi,\widetilde{\SL}_2)$ of $M$ is a trivial principal
${\mathbb Z}^{2g}$-bundle over the Weil space
$\R (\pi^{\orb},\PSL_2{\mathbb R})$ of~$O_M$.
\end{thm}

\begin{rem}
{\rm Our argument is largely analogous to considerations of
Oh\-shika~\cite{ohsh87} concerning the space $\R (\pi ,\Isom_0(
\widetilde{\SL}_2))$. However, he uses an invalid presentation for~$\pi$.
Our argument repairs that oversight.}
\end{rem}

We prepare for the proof of Theorem~\ref{thm:weil-sl} by first
establishing a concrete topological realisation of the orbifold
fundamental group, and by describing the geometry 
$\widetilde{\SL}_2$ in a little more detail.

\vspace{2mm}

\noindent {\bf The deck transformation group} $\pi^{\orb}$
\vspace{.5mm}

\noindent Let $O_M$ be a fixed topological orbifold of genus~$g$
and with $n$ cone points of multiplicity $\alpha_1,\ldots ,
\alpha_n$. We do not yet fix a hyperbolic structure on~$O_M$.
Choose a base point $x_0\in O_M$ distinct from all the cone points,
and a lift $\widetilde{x}_0\in\widetilde{O}_M$ of $x_0$
in the universal covering space~$\widetilde{O}_M$. Choose a system
of $2g$ loops on~$O_M$, based at~$x_0$, and a curve from $x_0$ to
each of the cone points, such that $O_M$ looks as in Figure~1
when cut open along these $2g+n$ curves. We may interpret that figure
as a fundamental region in~$\widetilde{O}_M$; it is determined
(amongst all possible fundamental regions whose boundary polygon
maps to the chosen system of curves) by the indicated placement of
$\widetilde{x}_0$ on its boundary. Notice that the sides
of this polygon labelled $\varepsilon_j,\varepsilon_j'$ meet
at a vertex mapping to the $j$-th cone point in~$O_M$; all other
vertices are lifts of~$x_0$.

\begin{figure}[h]
\centerline{\relabelbox\small
\epsfxsize 13cm \epsfbox{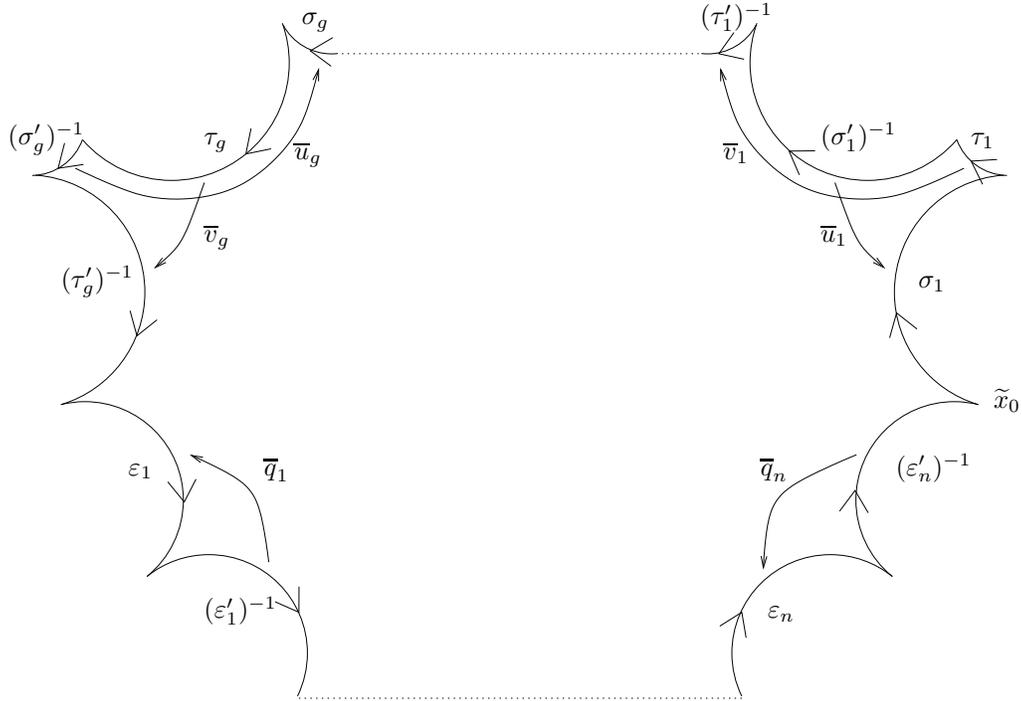}
\extralabel <-3.3cm, 3cm> {$\oq_n$}
\extralabel <-9.9cm, 3cm> {$\oq_1$}
\extralabel <-.2cm, 3.9cm> {$\widetilde{x}_0$}
\extralabel <-2.5cm, 6.1cm> {$\ou_1$}
\extralabel <-10.7cm, 6.1cm> {$\ov_g$}
\extralabel <-3.8cm, 7.2cm> {$\ov_1$}
\extralabel <-9.5cm, 7.2cm> {$\ou_g$}
\extralabel <-3.2cm, 1.1cm> {$\varepsilon_n$}
\extralabel <-10.7cm, 1.1cm> {$(\varepsilon_1')^{-1}$}
\extralabel <-1.5cm, 3cm> {$(\varepsilon_n')^{-1}$}
\extralabel <-11.7cm, 3cm> {$\varepsilon_1$}
\extralabel <-1.2cm, 5.5cm> {$\sigma_1$}
\extralabel <-12.6cm, 5.5cm> {$(\tau_g')^{-1}$}
\extralabel <-.5cm, 7.4cm> {$\tau_1$}
\extralabel <-13.3cm, 7.4cm> {$(\sigma_g')^{-1}$}
\extralabel <-2.5cm, 7.4cm> {$(\sigma_1')^{-1}$}
\extralabel <-10.7cm, 7.4cm> {$\tau_g$}
\extralabel <-4.1cm, 9cm> {$(\tau_1')^{-1}$}
\extralabel <-9.4cm, 9cm> {$\sigma_g$}
\endrelabelbox}
\caption{A fundamental domain for $O_M$.}
\end{figure}

Let $\ou_1,\ov_1,\ldots ,\ou_g,\ov_g,\oq_1,\ldots ,\oq_n$ be the
deck transformations of $\widetilde{O}_M$ which effect the gluing
maps of the sides of the chosen fundamental polygon as indicated
in Figure~1. The deck transformation
$\prod_i[\ou_i,\ov_i]\prod_j \oq_j$
(read from the right as a composition of maps) fixes the
point~$\widetilde{x}_0$, which is not the lift of a cone point,
so we conclude
\[ \prod_i[\ou_i,\ov_i]\prod_j \oq_j=1.\]
Similarly, we have
\[ \oq_j^{\alpha_j}=1,\;\; j=1,\ldots ,n.\]
So our choices provide us with a topological realisation of
$\pi^{\orb}$ as deck transformation group, consistent with
the presentation of $\pi^{\orb}$ fixed earlier on. Once $O_M$
is equipped with a hyperbolic structure and an orientation,
then $\widetilde{O}_M={\mathbb H}^2$ and the $\ou_i,\ov_i,\oq_j$
are isometries of~${\mathbb H}^2$, i.e.\ elements of $\PSL_2
{\mathbb R}$. The identification of $\widetilde{O}_M$ with
${\mathbb H}^2$ is uniquely determined if we specify, for instance,
the lift $\widetilde{x}_0\in{\mathbb H}^2$, the initial direction
of $\sigma_1$ at that point, and require that the orientation
lifted from $O_M$ coincides with a chosen orientation
of~${\mathbb H}^2$. In this way an oriented hyperbolic structure
on $O_M$ defines an element of $\R (\pi^{\orb},\PSL_2{\mathbb R})$.

\vspace{2mm}

\newpage
\noindent {\bf The geometry} $\widetilde{\SL}_2$
\vspace{.5mm}

\noindent
Identify $\widetilde{\SL}_2$, the universal cover of the
unit tangent bundle $\ST {\mathbb H}^2$, with ${\mathbb H}^2\times {\mathbb R}$
with coordinates $(z,\theta)$, cf.~\cite{wall86}, \cite{ue91}, \cite{gego95}.
An element $\overline{A}
=\left(\begin{array}{cc}a&b\\c&d\end{array}\right)\in\PSL_2{\mathbb R}$
acts isometrically on ${\mathbb H}^2$ by the fractional linear transformation
\[ \overline{A}(z)=\frac{az+b}{cz+d}.\]
On differentiation this yields
\[ \partial_z\longmapsto \frac{\partial_z}{(cz+d)^2}.\]
By taking the argument modulo $2\pi$ of the logarithm of this action,
we obtain the description of the differential of $\overline{A}$ in terms of
the coordinates $(z,\theta\;\mbox{\rm mod}\; 2\pi)$ for $\ST {\mathbb H}^2$:
\[ (z,\theta\;\mbox{\rm mod}\; 2\pi )\longmapsto
\left( \frac{az+b}{cz+d},\theta-2\arg (cz+d)\;
\mbox{\rm mod}\; 2\pi\right).\]
An element $A\in\widetilde{\SL}_2$ that maps to $\overline{A}\in\PSL_2
{\mathbb R}$ then acts by left-multiplication on $\widetilde{\SL}_2=
\widetilde{\ST {\mathbb H}^2}$ as
\[ A(z,\theta )=\left( \frac{az+b}{cz+d},\theta-2\arg (cz+d)\right) ,\]
where the argument function depends on $z\in{\mathbb H}^2$ and the
lift $A$ of $\overline{A}$.

For $\overline{A}$ hyperbolic, or elliptic of finite order, we can define
a {\em preferred lift} $A$ as follows:

\begin{itemize}
\item If $\overline{A}$ is a hyperbolic element, then the action of
$A$ on the $\theta$-coordinate is defined by parallel translation of
a unit tangent vector along the axis of~$\overline{A}$.

\item If $\overline{A}$ is an elliptic element of
order~$\alpha\in {\mathbb N}$,
rotating by $2\pi/\alpha$ in positive direction around its fixed point
(with respect to the
complex orientation of~${\mathbb H}^2$), then $A$ is determined by
$A^{\alpha}(z,\theta )=(z,\theta +2\pi )$.
\end{itemize}

More geometrically, this means the following. Given a hyperbolic
element~$\overline{A}$, take the path of hyperbolic
elements, all with the same axis, from the identity in $\PSL_2{\mathbb R}$
to $\overline{A}$. If $\overline{A}$ is elliptic, rotating by $2\pi/\alpha$
in positive direction, take the path of rotations about the same fixed point,
starting at the identity and ending with~$\overline{A}$ (without completing
any full turns). The lifts of these paths to~$\widetilde{\SL}_2$,
starting at the identity, end at the preferred lift $A$
of~$\overline{A}$.

The sign in the lift
of an elliptic element is crucial for the subsequent discussion, so we give
a brief justification for it. Consider the element
\[ \overline{A}_t=\left(\begin{array}{rr}\cos t&\sin t\\ -\sin t&\cos t
\end{array}\right) .\]
Notice that, for increasing $t\in {\mathbb R}$, the action of
$\overline{A}_t$ on ${\mathbb H}^2$ by fractional linear transformations
defines a positive (i.e.\ counterclockwise)
rotation about $i\in{\mathbb H}^2$.
The differential of $\overline{A}_t$ is given by
\[ (z,\theta )\longmapsto (\overline{A}_t(z),\theta-2\arg (-z\sin t+\cos t)
\;\mbox{\rm mod}\; 2\pi ). \]
In particular,
\[ (i,\theta )\longmapsto (i,\theta -2\arg (-i\sin t+\cos t)\;
\mbox{\rm mod}\; 2\pi ). \]
Up to integer multiples of $2\pi$ we have $\arg (\cos t -i\sin t)=-t$.
Therefore the choice of $-t$ for this argument is the only one
for which the lifts $A_t$ of $\overline{A}_t$ form a $1$-parameter
subgroup. This yields
\[ A_{\pi}(z,0)=(z,\theta +2\pi ).\]

Recall that the identity component $\Isom_0(\widetilde{\SL}_2)$ of the full
isometry group of $\widetilde{\SL}_2$ is $\widetilde{\SL}_2\times_{\mathbb Z}
{\mathbb R}$, where the ${\mathbb R}$-factor acts by translation of the
$\theta$-coordinate, and ${\mathbb Z}$ corresponds to the kernel of the
$(\widetilde{\SL}_2\times {\mathbb R})$-action, generated by
$(A_{\pi},-2\pi )$, cf.~\cite[p.~188]{gego95} (beware that in the cited
reference we worked with a negative rotation).

\vspace{2mm}

\noindent {\it Proof of Theorem}~\ref{thm:weil-sl}.
Suppose we are given $\orh\in\R (\pi^{\orb},
\PSL_2{\mathbb R})$. Assume that the orientation of $O_M$
is chosen in such a way
that the loop
\begin{equation}
\label{eqn:loop}
\prod_i\sigma_i\tau_i(\sigma_i')^{-1}(\tau_i')^{-1}\prod_j
\varepsilon_j(\varepsilon_j')^{-1}
\end{equation}
(where the composition of paths is read, as usual, from
the left),
corresponds under the developing
map to a positively oriented polygonal loop in~${\mathbb H}^2$. This amounts
to fixing one of the two components of $\R (\pi^{\orb},\PSL_2{\mathbb R})$,
and geometrically it means that Figure~1 may be interpreted as
a fundamental region inside ${\mathbb H}^2$ with its standard orientation.
In particular, the $\oq_j$ (or more precisely: $\orh (\oq_j)$) are
rotations in positive direction by $2\pi /\alpha_j$.

Now lift $\orh$ to a representation $\rho_0\co
\pi\rightarrow\Isom_0(\widetilde{\SL}_2)$ as follows:

\begin{itemize}
\item The $\rho_0(u_i)$ and $\rho_0(v_i)$ are defined as the preferred
lifts of $\orh (\ou_i)$ and $\orh (\ov_i)$, respectively.
\item $\rho_0(h)(z,\theta )=(z,\theta +2\pi r)$.
\item $\rho_0(q_j)$ is the lift of $\orh (\oq_j)$ that satisfies
$\rho_0(q_j)^{\alpha_j}\rho_0(h)^{\beta_j}=1$.
\end{itemize}
(For the other choice of orientation, the sign in the
definition of $\rho_0(h)$ needs to be changed).

\begin{rem}
{\rm The integer $r$ is what we called the {\em fibre index} in~\cite{gego95a}.
Geo\-metrically it describes the number of full turns along
a regular fibre of $M\rightarrow O_M$ made by the two contact planes
in the standard taut contact circle on
$\rho_0(\pi)\backslash\widetilde{\SL}_2$.}
\end{rem}

\begin{lem}
\label{lem:rep}
This does indeed define a representation of~$\pi$, and in fact 
we have $\rho_0\in\R (\pi ,\widetilde{\SL}_2)$.
\end{lem}

\noindent {\em Proof.}
We first check that $\rho_0(q_j)\in\widetilde{\SL}_2$; for all the other
generators of $\pi$ this is clear by definition. The preferred lift of
$\orh (\oq_j)$ lies in $\widetilde{\SL}_2$. With our choice of
orientation, $\orh (\oq_j)$ is an elliptic element rotating by
$2\pi /\alpha_j$ in positive direction around its fixed point. So the
$\alpha_j$-th power of the preferred lift of $\orh (\oq_j)$ gives
a shift by $2\pi$ in $\theta$-direction.
Thus $\rho_0(q_j)$ differs
from the preferred lift by a translation in $\theta$-direction by
\[ -2\pi \frac{1+\beta_jr}{\alpha_j}=-2\pi (1+k_j)\in
2\pi {\mathbb Z},\]
so it also lies in~$\widetilde{\SL}_2$. 

Provided that $\rho_0$ is a representation of $\pi$, it is clear from the
geometry of $\widetilde{\SL}_2$ that $\rho_0(\pi)\backslash\widetilde
{\SL}_2$ will define a Seifert manifold $M$ with base orbifold~$O_M$.
So to verify that $\rho_0\in\R (\pi,\widetilde{\SL}_2)$, it only remains
to be checked that our definition of $\rho_0$ is compatible with
the relations in~$\pi$.

Write $\rho_0'(q_j)$ for the lift of $\orh (\oq_j)$ to $\Isom_0
(\widetilde{\SL}_2)$ which satisfies $\rho_0'(q_j)^{\alpha_j}=1$
(that is, $\rho_0'(q_j)$ is distinguished amongst all lifts
of $\orh (\oq_j)$ as the one which does {\em not} act as a non-trivial
helicoidal motion).
This $\rho_0'(q_j)$ differs from the preferred lift of $\orh (\oq_j)$ by
a shift in $\theta$-direction by~$-2\pi /\alpha_j$; in particular,
$\rho_0'(q_j)\not\in\widetilde{\SL}_2\subset\Isom_0(\widetilde{\SL}_2)$.
The reason for introducing this lift is its geometric significance
implicit in the following statement:

\begin{equation}
\label{eqn:parallel}
\prod_i[\rho_0(u_i),\rho_0(v_i)]
\prod_j\rho_0'(q_j)(z,\theta )=(z,\theta-2\pi\chi^{\orb}(O_M)).
\end{equation}

This equation
also appears in~\cite{ue91}, but we trust the reader will appreciate
a little more elucidation than is given there.
Notice that $-2\pi\chi^{\orb}$ is the area of~$O_M$, hence may
be thought of as {\em minus}
the integral of the constant Gau{\ss} curvature~$-1$
over a fundamental region for $O_M$ in~${\mathbb H}^2$. To prove
equation~(\ref{eqn:parallel})
it therefore suffices to show, by the classical Gau{\ss}-Bonnet
theorem, that the total effect on the $\theta$-coordinate
of the isometries on the left-hand side is equal to the parallel
transport along the inverse of the loop~(\ref{eqn:loop}).

Different choices of lifts of $\orh (\ou_i)$ and $\orh (\ov_i)$ do not
affect the commutators on the left-hand side of~(\ref{eqn:parallel}),
so instead of $\rho_0(u_i)$ and $\rho_0(v_i)$ we may use lifts
$\rho_0'(u_i)$ and $\rho_0'(v_i)$ defined as follows (refer to
Figure~2 for notation, where we drop the index~$i$ for convenience,
and recall that $\widetilde{\SL}_2=\widetilde{\ST {\mathbb H}^2}$).

\begin{itemize}
\item For $X\in\widetilde{\ST_c{\mathbb H}^2}$, that is, an element
of $\widetilde{\ST {\mathbb H}^2}$ in the ${\mathbb R}$-fibre
over $c\in{\mathbb H}^2$, we have that $\rho_0'(u)(X)\in
\widetilde{\ST_b{\mathbb H}^2}$ is obtained by parallel
transport of $X$ along $\tau^{-1}$, and $\rho_0'(v)(X)\in
\widetilde{\ST_d{\mathbb H}^2}$ is obtained by parallel transport of $X$
along~$(\sigma ')^{-1}$.
\end{itemize}

\begin{figure}[h]
\centerline{\relabelbox\small
\epsfxsize 9cm \epsfbox{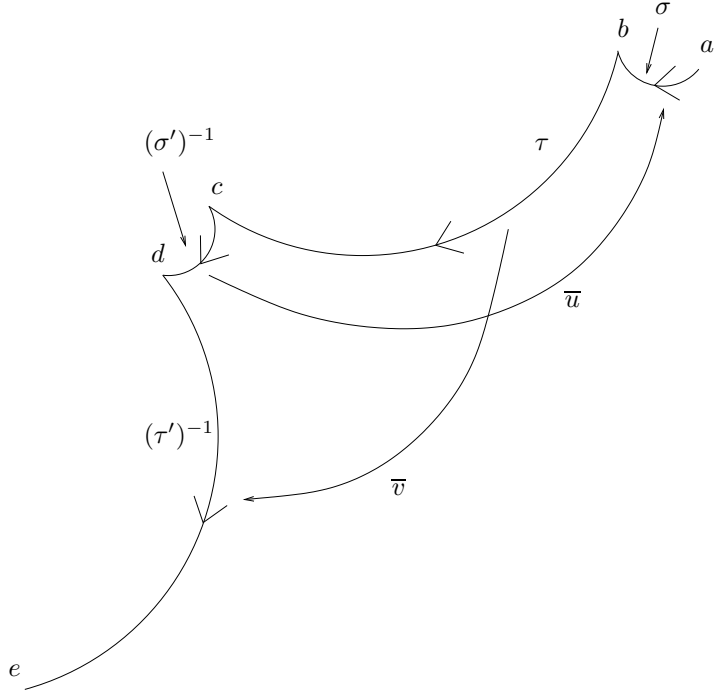}
\extralabel <0cm, 8.5cm> {$a$}
\extralabel <-1.1cm, 8.7cm> {$b$}
\extralabel <-6.5cm, 6.6cm> {$c$}
\extralabel <-7.3cm, 5.7cm> {$d$}
\extralabel <-9.2cm, .2cm> {$e$}
\extralabel <-1.8cm, 5.1cm> {$\ou$}
\extralabel <-4.1cm, 2.6cm> {$\ov$}
\extralabel <-7.4cm, 3.3cm> {$(\tau ')^{-1}$}
\extralabel <-2.2cm, 7.2cm> {$\tau$}
\extralabel <-.6cm, 9cm> {$\sigma$}
\extralabel <-7.4cm, 7.2cm> {$(\sigma ')^{-1}$}
\endrelabelbox}
\caption{Definition of $\rho_0'(u)$ and $\rho_0'(v)$.}
\end{figure}

Since $\orh (\ou )$ maps $\sigma '$ isometrically to $\sigma$ and
$\orh (\ov )$ maps $\tau$ isometrically to~$\tau '$, it is
a simple matter to check that the map
\[ [\rho_0'(u),\rho_0'(v)]\co \widetilde{\ST_e{\mathbb H}^2}
\longrightarrow \widetilde{\ST_a{\mathbb H}^2}\]
is given by parallel transport along the path
$\tau '\sigma '\tau^{-1}\sigma^{-1}$.

For the elliptic elements we argue similarly. The previously defined
$\rho_0'(q_j)$ fixes the fibre of $\widetilde{\ST {\mathbb H}^2}$
over the corresponding cone point, and since $\orh (\oq_j)$
maps $\varepsilon_j'$ isometrically onto~$\varepsilon_j$, we see
that $\rho_0'(q_j)$ is given by parallel transport along~$\varepsilon_j'
\varepsilon_j^{-1}$. This proves equation~(\ref{eqn:parallel}).

We can now complete the proof of Lemma~\ref{lem:rep}. We only need to verify
that $\rho_0$ respects the relation
$\prod_i[u_i,v_i]\prod_jq_j=h^b$.
We have
\[ \rho_0(q_j)(z,\theta )=\rho_0'(q_j)(z,\theta )-(0,2\pi
\frac{\beta_jr}{\alpha_j}).\]
Hence
\[ \prod_i[\rho_0(u_i),\rho_0(v_i)]\prod_j\rho_0(q_j)(z,\theta)=
(z,\theta+\theta_0)\]
with
\begin{eqnarray*}
\theta_0/2\pi & = & -\chi^{\orb}(O_M)-\sum_j\frac{\beta_jr}{\alpha_j}\\
     & = & -\Bigl( 2-2g-n+\sum_j\frac{1}{\alpha_j}\Bigr) -
             \sum_j\Bigl(1-\frac{1}{\alpha_j}+k_j\Bigr) \\
     & = & 2g-2-\sum_j k_j \\
     & = & rb,
\end{eqnarray*}
thus
\[ \prod_i[\rho_0(u_i),\rho_0(v_i)]\prod_j\rho_0(q_j)=\rho_0(h)^b,\]
as desired. This completes the proof of Lemma~\ref{lem:rep}.\hfill $\Box$

\vspace{2mm}

We have thus found a particular lift $\rho_0\in\R (\pi ,\widetilde{\SL}_2)$
of $\orh\in\R (\pi^{\orb},\PSL_2{\mathbb R})$. It remains to settle
the question what ambiguity there is in the choice of such a lift. The
next lemma shows that for the generators $h,q_1,\ldots ,q_n$ there
is no ambiguity at all.

\begin{lem}
\label{lem:uniquelift}
Let $\rho\in\R (\pi,\widetilde{\SL}_2)$ be any other lift of~$\orh$.
Then $\rho (h)=\rho_0(h)$ and $\rho (q_j)=\rho_0(q_j)$.
\end{lem}

\noindent {\em Proof.}
Write
\begin{eqnarray*}
\rho (h)(z,\theta ) & = & \rho_0(h)(z,\theta)+(0,2\pi y),\\
\rho (q_j)(z,\theta ) & = & \rho_0(q_j)(z,\theta )+ (0,2\pi x_j).
\end{eqnarray*}
The relations $q_j^{\alpha_j}h^{\beta_j}=1$ and $\prod_i[u_i,v_i]
\prod_jq_j=h^b$ then yield the equation
\[ C\left(\begin{array}{c}x_1\\ \vdots\\x_n\\y\end{array}\right) =0\]
with
\[ C=\left(\begin{array}{ccccc}
\alpha_1 & 0 & \cdots & 0 &\beta_1\\
 0 & \alpha_2 & \cdots & 0 & \beta_2\\
 \vdots & \vdots & \ddots  & \vdots & \vdots\\
 0 & 0 & \cdots & \alpha_n & \beta_n\\
 1 & 1 & \cdots & 1 & -b
\end{array}\right) .\]
Recall that the Euler number $e$ was defined as $e=-(b+
\sum_j\beta_j/\alpha_j)$, hence
\begin{eqnarray*}
\det C & = & \left|
\begin{array}{ccccc}
\alpha_1 & 0 & \cdots & 0 &\beta_1\\
 0 & \alpha_2 & \cdots & 0 & \beta_2\\
 \vdots & \vdots & \ddots  & \vdots & \vdots\\
 0 & 0 & \cdots & \alpha_n & \beta_n\\
 0 & 0 & \cdots & 0 & e
\end{array}\right| \\
       & = & \alpha_1...\alpha_ne\\
       & \neq & 0.
\end{eqnarray*}
This implies that the only solution is $x_1=...=x_n=y=0$.\hfill $\Box$

\begin{rem}
{\rm Observe that the condition $e\neq 0$ (that is,
$\widetilde{\SL}_2$-geometry as opposed to $M$ being modelled
on ${\mathbb H}^2\times{\mathbb R}$) enters crucially in this lemma.
Notice further that $\rho (h)$ is forced to be the vertical
shift by~$2\pi r$, with $r$ equal to the fibre index, no matter
what $\orh$ is.}
\end{rem}

For the generators $u_i,v_i$, on the other hand, there is a certain
freedom in choosing lifts. Indeed,
we can define a lift $\rho_{\mathbf w}\in
\R (\pi ,\widetilde{\SL}_2)$ of the representation
$\orh\in\R (\pi^{\orb},\PSL_2{\mathbb R})$
for any ${\mathbf w}\in{\mathbb Z}^{2g}$ by setting
\begin{eqnarray*}
\rho_{\mathbf w}(u_i)(z,\theta ) & = & \rho_0(u_i)(z,\theta )+
                                      (0,2\pi w_{2i-1}),\\
\rho_{\mathbf w}(v_i)(z,\theta ) & = & \rho_0(v_i)(z,\theta )+
                                      (0,2\pi w_{2i}),
\end{eqnarray*}
and any lift of $\orh$ to $\R (\pi ,\widetilde{\SL}_2)$ must be of this form.

Conversely, if we are given a representation $\rho\in\R (\pi ,
\widetilde{\SL}_2)$, the discrete subgroup $\rho (\pi )\subset
\widetilde{\SL}_2$ preserves the fibration $\widetilde{\SL}_2
\rightarrow {\mathbb H}^2$, which induces a Seifert fibration of
$\rho (\pi )\backslash \widetilde{\SL}_2$. The central element $h\in\pi$
must map to the centre of~$\widetilde{\SL}_2$, which equals the kernel
of the projection
\[ \Isom_0(\widetilde{\SL}_2)\supset\widetilde{\SL}_2
\stackrel{\mbox{\scriptsize\rm pr}}{\longrightarrow}\Isom_0 {\mathbb H}^2
=\PSL_2{\mathbb R},\]
and $\mbox{\rm pr}\circ\rho$ defines an element of $\R (\pi^{\orb},
\PSL_2{\mathbb R})$, cf.~\cite[Thm.~4.15]{scot83}.

Notice that we have a free ${\mathbb Z}^{2g}$-action on $\R (\pi ,
\widetilde{\SL}_2)$ defined by
\[ \begin{array}{ccc}
{\mathbb Z}^{2g}  \times  \R (\pi ,\widetilde{\SL}_2) &
\longrightarrow & \R (\pi ,\widetilde{\SL}_2)\\
( {\mathbf w'},\rho_{\mathbf w}) & \longmapsto & \rho_{{\mathbf w}
+{\mathbf w'}},
\end{array} \]
and a section $\rho_0$ of the covering $\R(\pi ,\widetilde{\SL}_2)
\rightarrow\R (\pi^{\orb},\PSL_2{\mathbb R})$.

This concludes the
proof of Theorem~\ref{thm:weil-sl}.\hfill $\Box$

\begin{rem}
\label{rem:r}
{\rm Observe that with $C$ the matrix from the proof of
Lemma~\ref{lem:uniquelift}, $|\det C|=\alpha_1\cdots \alpha_n|e|$
equals the order of the torsion subgroup of the abelia\-nised~$\pi$.
The $\alpha_j$ are recovered from $\pi$ as the orders of the different
maximal cyclic subgroups (modulo conjugation) of~$\pi^{\orb}$, which
is the quotient of $\pi$ by its centre $\langle h\rangle$. The genus
$g$ of the base orbifold is determined by the rank of the
abelia\-nised~$\pi$ (which equals~$2g$). This allows to compute
$|r|=|\chi^{\orb}/e|$ from the abstract group structure of $\pi$
alone. Fixing a sign of $r$ amounts to fixing orientations of
$O_M$ and the Seifert fibres; with this choice of sign the
normalised Seifert invariants are recovered from~$\pi$. Alternatively,
the Seifert invariants and $\pi$ can be recovered from
$\pi^{\orb}$ and~$r$.}

\end{rem}

\begin{rem}
{\rm Throughout this paper we have been thinking of $\pi$ and
$\pi^{\orb}$ as deck transformation groups, rather than as
fundamental groups of homotopy classes of based loops. In the
course of proving equation~(\ref{eqn:parallel}), however, we had
to show that the action (on the fibre of $\widetilde{\SL}_2
\rightarrow {\mathbb H}$) of a certain word $w\in\pi$ could
be computed as parallel transport along a suitable loop
in~${\mathbb H}$. It is tempting to believe that this loop
corresponds directly to the projected word $\overline{w}\in
\pi^{\orb}$ if $\pi^{\orb}$ is reinterpreted in terms of homotopy
classes of loops. We deem it worth pointing out that the relation
between the fundamental group of a topological space $X$
interpreted as a deck transformation group (which we shall
denote~$\pi$) and the fundamental group interpreted in terms
of loops based at~$x_0$ (for which we write $\pi_1(X,x_0)$) is
a little more subtle, a fact that is usually brushed over in
elementary treatments. (For orbifold fundamental groups there
are completely analogous statements.)

Write $p\co\widetilde{X}\rightarrow X$ for the universal covering.
It is well-known that $\pi_1(X,x_0)$ acts on $p^{-1}(x_0)$
from the {\em right}, with $[\gamma ]\in\pi_1(X,x_0)$ sending
$y\in p^{-1}(x_0)$ to the endpoint $y\gamma$ of the lift of
$\gamma$ to $\widetilde{X}$ with initial point~$y$. If we fix an
$\widetilde{x}_0\in p^{-1}(x_0)$, there is a unique deck transformation
$u_{\gamma}$ such that $u_{\gamma}(\widetilde{x}_0)=\widetilde{x}_0
\gamma$. One checks easily that the map $[\gamma ]\mapsto u_{\gamma}$
defines an isomorphism $\pi_1(X,x_0)\rightarrow\pi$, where
composition in $\pi_1(X,x_0)$ is read from the left as composition
of paths, in $\pi$ from the right as composition of maps.

However, this isomorphism depends crucially on the choice
of~$\widetilde{x}_0$. If $\widetilde{x}_0$ is replaced by
$u(\widetilde{x}_0)$ for some $u\in\pi$, one verifies that the
corresponding isomorphism $\pi_1(X,x_0)\rightarrow \pi$ is given
by $[\gamma ]\mapsto uu_{\gamma}u^{-1}$.

This implies that the loop~(\ref{eqn:loop}) corresponds, under the
isomorphism of $\pi_1^{\orb}(O_M,x_0)$ with $\pi^{\orb}$ determined
by the choice of~$\widetilde{x}_0$, not to the word in
equation~(\ref{eqn:parallel}), but to a word where the letters
have been replaced by increasingly complicated conjugate letters.}
\end{rem}

\begin{thm}
\label{thm:deform-sl}
Let $M$ be a left-quotient of $\widetilde{\SL}_2$. The Teichm\"uller
space $\T (M)$ of taut contact circles on~$M$ is a trivial principal
${\mathbb Z}^{2g}$-bundle over the Teich\-m\"uller space $\T (O_M)$ of
hyperbolic metrics on the base orbifold $O_M$ of the unique Seifert fibration
$M\rightarrow O_M$.

The moduli space $\M (M)$ is an $r^{2g}$-fold branched covering of~$\M(O_M)$,
where $g$ is the genus of $O_M$ and
the integer $r$ is determined from the Seifert bundle structure
$M\rightarrow O_M$ by the equation $r=\chi^{\orb}(O_M)/e$.
\end{thm}

\noindent {\em Proof.}
We have seen that the projection $\widetilde{\SL}_2\rightarrow
\PSL_2{\mathbb R}$ induces a covering map $\R (\pi ,\widetilde{\SL}_2)
\rightarrow \R (\pi^{\orb},\PSL_2{\mathbb R})$, and this in turn
induces a map $\T^{\alg}(M)\rightarrow\T^{\alg}(O_M)$. This map is
well-defined, since any inner automorphism of $\widetilde{\SL}_2$
induces an inner automorphism of $\PSL_2{\mathbb R}$. Furthermore,
integer shifts in fibre direction of $\widetilde{\SL}_2\rightarrow
{\mathbb H}^2$ lie in the centre of $\widetilde{\SL}_2$, so it follows
directly from the proof of Theorem~\ref{thm:weil-sl} that $\T^{\alg}(M)$
is a trivial principal ${\mathbb Z}^{2g}$-bundle over $\T (O_M)$.

Similarly, we have a map $\T (M)\rightarrow \T(O_M)$, defined as follows:
The common kernel $\ker\omega_1\cap\ker\omega_2$ of a taut contact circle
gives $M$ the structure of a Seifert fibration. The Seifert fibration
of $M$ being unique up to isotopy, we can choose a representative of
$[(\omega_1,\omega_2)]\in\T (M)$ such that it induces the fixed Seifert
fibration $M\rightarrow O_M$. The unique $(-1)$-Cartan structure in the
homothety class of this representative then induces a hyperbolic metric
and an orientation on~$O_M$.

It is not difficult to see that these maps fit together to form a commutative
diagram
\[
\begin{CD}
\T (M) @>\Phi>> \T^{\alg}(M)\\
@VVV             @VVV\\
\T (O_M) @>\phi>> \T^{\alg}(O_M).
\end{CD}
\]

\begin{rem}
{\rm There are natural topologies on $\T (M)$ and $\T^{\alg}(M)$ which turn
the vertical arrows into covering maps. The fact that $\phi$ is a homeomorphism
then implies that the bijection $\Phi$ is a homeomorphism as well.}
\end{rem}

In order to determine the moduli space $\M (M)$ we need some information
about the automorphism group of~$\pi$.

\begin{lem}
\label{lem:aut}
There is a split short exact sequence
\[ 0\longrightarrow {\mathbb Z}^{2g}\longrightarrow \Aut (\pi )\longrightarrow
\Aut(\pi^{\orb})\longrightarrow 1,\]
where ${\mathbf c}\in{\mathbb Z}^{2g}$ acts on $\pi$ by
\begin{eqnarray*}
u_i & \longmapsto & u_ih^{c_{2i-1}},\\
v_i & \longmapsto & v_ih^{c_{2i}}.
\end{eqnarray*}
\end{lem}

\noindent {\em Proof.}
Except for the surjectivity of the homomorphism $\Aut (\pi )\rightarrow
\Aut (\pi^{\orb})$ and the existence of a splitting
this follows easily from the explicit
presentations of $\pi$ and~$\pi^{\orb}$. This surjectivity is a property
specific to left-quotients of~$\widetilde{\SL}_2$.

We argue geometrically, as in the proof of the surjectivity of $\Phi$
in Section~\ref{section:Cartan}. Given an element of
$\Aut (\pi^{\orb})$, it can be realised geometrically by a diffeomorphism
$\overline{\psi}$ of~$O_M$, which may only permute cone points of
equal multiplicity. We may assume that $\overline{\psi}$ fixes
an additional base point, distinct from any of the cone points.
Changing, if necessary, $\overline{\psi}$ by an isotopy fixing the base
point $x_0$ and the cone points $x_1,\ldots ,x_n$, we can find disjoint
disc neighbourhoods $D_0,\ldots ,D_n$ of these points such that
$\overline{\psi}$ preserves the complement $U\subset O_M$ of this
collection of discs and either fixes the boundary curve of each
disc (if $\overline{\psi}$ is orientation preserving) or reverses
each of these curves (in the orientation reversing case). Over
$U$ the Seifert bundle is trivial, so $\overline{\psi}$ lifts to
the diffeomorphism $\overline{\psi}\times (\pm\mbox{\rm id})$ of
$U\times S^1$.

By Remark~\ref{rem:mult} the Seifert invariant $(\alpha ,\beta )$ of
a multiple fibre only depends on the multiplicity~$\alpha$ (for a fixed $M$
with a choice of base and fibre orientations,
and hence fixed~$r$). Depending on whether $\overline{\psi}$ preserves
or reverses orientation, it is a straightforward matter to check that
exactly one of $\overline{\psi}\times (\pm\mbox{\rm id})$ extends
to a diffeomorphism $\psi$ of~$M$
which induces an automorphism of $\pi$ that projects to the given
automorphism of~$\pi^{\orb}$.

For a fixed choice of trivialisation $U\times S^1$, the diffeomorphism
$\psi$ is determined up to isotopy by $\overline{\psi}$, so it follows
that the map $\Aut (\pi^{\orb})\rightarrow\Aut (\pi )$, associating
to the automorphism induced by $\overline{\psi}$ the one induced
by~$\psi$, is a splitting of the short exact sequence. Since $h$ is
central in~$\pi$, one sees for purely algebraic reasons that there
is a ${\mathbb Z}^{2g}$-worth of such splittings. Geometrically,
they correspond to the different choices of
trivialisations $U\times S^1$. \hfill $\Box$

\vspace{2mm}

We continue with the proof of Theorem~\ref{thm:deform-sl}.
Before turning to the moduli space, we consider the Chabauty space
\[ {\mathcal S}(M)=\R (\pi ,\widetilde{\SL}_2)/\Aut (\pi ). \]
By Lemma~\ref{lem:uniquelift} we have $\rho (h)(z,\theta )=
(z,\theta +2\pi r)$ for any $\rho\in\R (\pi ,\widetilde{\SL}_2)$, and from
the definition of the principal ${\mathbb Z}^{2g}$-action on
$\R (\pi ,\widetilde{\SL}_2)$ we infer that ${\mathbf c}\in
{\mathbb Z}^{2g}\subset \Aut (\pi )$ acts on $\R (\pi ,
\widetilde{\SL}_2)$ by
\[ \rho_{\mathbf w}\longmapsto \rho_{{\mathbf w}+r{\mathbf c}}.\]
Hence, by first dividing out this ${\mathbb Z}^{2g}$-action, and then
the $\Aut(\pi^{\orb})$-action, we get
\[ {\mathcal S}(M)=(\R (\pi^{\orb},\PSL_2{\mathbb R})\times 
{\mathbb Z}_r^{2g})/\Aut(\pi^{\orb}).\]
The Weil space consists of faithful representations, and so
$\Aut (\pi^{\orb})$ acts freely on it. It follows that ${\mathcal S}(M)$
is a ${\mathbb Z}_r^{2g}$-bundle, i.e.\ an $r^{2g}$-fold covering, of
\[ {\mathcal S}(O_M)=(\R (\pi^{\orb},\PSL_2{\mathbb R}))/\Aut(\pi^{\orb}).\]

The moduli space $\M (M)$ may now be regarded as 
\[ \T^{\alg}(M)/\Out (\pi )=(\T^{\alg}(M)/{\mathbb Z}^{2g})/
\Out(\pi^{\orb}),\]
or as
\[  \Inn (\widetilde{\SL}_2)\backslash {\mathcal S}(M).\]
In either case, $\M (M)$ is obtained as the quotient under an action
that may stabilise a point in the base (by which we mean $\T^{\alg}(O_M)$
or ${\mathcal S}(O_M)$, respectively), without stabilising the points
sitting over it, cf.~\cite[pp.~194--195]{klr85}. This results in the
branching of the covering $\M (M)$ over $\M(O_M)$.

Geometrically, points with non-trivial
stabiliser correspond to Seifert manifolds resp.\ base orbifolds
with non-trivial isometries. \hfill $\Box$

\vspace{2mm}

A fundamental question about the moduli space for any kind of geometric
structure is whether it is connected. The moduli space being connected
amounts to there being a unique structure up to diffeomorphism and
deformations.

In our situation, this issue can be addressed by studying the covering
$\T^{\alg}(M)/{\mathbb Z}^{2g}\rightarrow\T (O_M)$; details will
appear in~\cite{gegob}. The fibre of this covering can be shown to
be in natural one-to-one correspondence with the fibrewise $r$-fold
coverings $M\rightarrow STO_M$, with $STO_M$ denoting the unit
tangent bundle of~$O_M$ (which is defined even in the presence of
cone points). In particular, for $r=2$ and $O_M$ a surface without
cone points, this corresponds to spin structures on~$O_M$. More generally,
for a choice of a hyperbolic structure on $O_M$ and an $r$-fold
fibrewise covering $M\rightarrow STO_M$, sections of $M\rightarrow O_M$
are roots of index $r$ of (unit) tangent vectors, just like
automorphic forms of weight $1/r$ for the subgroup
of $\PSL_2{\mathbb R}$ defining~$O_M$. The latter are traditionally
required to correspond to holomorphic sections
of the associated complex line
bundle; in the presence of cone points, this last statement needs
to be qualified.

In the case that $O_M$ is free of cone points, it turns out that
$\M (M)$ is connected for $r$ odd, and it has exactly two components
for $r$ even. For $r=2$ this corresponds to the fact that there
are exactly two spin structures on $O_M$ up to diffeomorphism,
cf.~\cite{dape86}. It is worth pointing out that any surface of genus
$g\geq 2$ (in fact, $g\geq 1$ suffices) with its two non-diffeomorphic
spin structures represents the two elements in the spin cobordism
group $\Omega^{\Spin}_2\cong{\mathbb Z}_2$.

We conclude that there is indeed non-trivial branching in
the covering $\M (M)\rightarrow \M (O_M)$.
\section{$\widetilde{\E}_2$-geometry}
\label{section:e}
The left-quotients of $\widetilde{\E}_2$ under a discrete, cocompact
subgroup $\Gamma$ are, up to diffeomorphism, exactly the five $T^2$-bundles
over $S^1$ with periodic monodromy $A\in\SL_2{\mathbb Z}$, cf.~\cite{rava81}.
Table~1 gives a choice of possible monodromy matrices~$A$, their periods,
and the lattices in ${\mathbb C}$ generated by $1$ and $\tau$,
$\mbox{\rm Im}\, \tau >0$, which are invariant under the action of
$A\in\SL_2{\mathbb Z}\subset\GL_2{\mathbb R}$ (here the action is the usual
linear one on ${\mathbb R}^2={\mathbb C}$).

We label the five left-quotients of $\widetilde{\E}_2$ as $M_k$, with
$k$ denoting the period of the monodromy matrix
defining~$M_k$, but we continue to write $T^3$ for the $3$-torus~$M_1$.

If the monodromy matrix is $A=\left(\begin{array}{cc}\alpha & \beta\\
\gamma &\delta\end{array}\right)$, then the fundamental group $\pi =\pi_1(M)$
of the total space has the presentation
\[ \pi =\left\{ s,t,b\co st=ts,\; bsb^{-1}=s^{\alpha}t^{\gamma},\;
btb^{-1}=s^{\beta}t^{\delta}\right\} ,\]
cf.~\cite{safu83}. In Table~1 we also list the first homology
groups~$H_1(M)$, obtained by abelianising~$\pi_1(M)$. This gives a proof
of our claim made in Section~\ref{section:Cartan} that these spaces
are distinguished by their fundamental groups.

\begin{table}
\centering
\begin{tabular}{c|c|c|c}
$A$ & period & $\tau$ & $H_1(M)$\\ \hline
$\left(\begin{array}{rr}1&0\\0&1\end{array}\right)$ & $1$ & any &
    ${\mathbb Z}^3$ \rule{0cm}{1cm}\\
$\left(\begin{array}{rr}-1&0\\0&-1\end{array}\right)$ & $2$ & any &
    ${\mathbb Z}\oplus{\mathbb Z}_2\oplus{\mathbb Z}_2$ \rule{0cm}{1cm}\\
$\left(\begin{array}{rr}0&-1\\1&-1\end{array}\right)$ & $3$ & $\exp (2\pi i/3)$
& ${\mathbb Z}\oplus {\mathbb Z}_3$ \rule{0cm}{1cm}\\
$\left(\begin{array}{rr}0&-1\\1&0\end{array}\right)$ & $4$ & $i$ &
    ${\mathbb Z}\oplus{\mathbb Z}_2$ \rule{0cm}{1cm}\\
$\left(\begin{array}{rr}0&-1\\1&1\end{array}\right)$ & $6$ & $\exp (2\pi i/6)$
 & ${\mathbb Z}$ \rule{0cm}{1cm}\\
\end{tabular}
\caption{$T^2$-bundles over $S^1$.}
\end{table}

We now want to study the Weil spaces $\R (\pi,\widetilde{\E}_2)$ for the five
different groups~$\pi$. Given a representation $\rho\co\pi\rightarrow
\widetilde{\E}_2$, write
\[ \overline{\rho}\co\pi\longrightarrow\E_2= \Isom_0{\mathbb R}^2 \]
for the composition of $\rho$ with the
projection $\widetilde{\E}_2\rightarrow\E_2$.

Recall the definition of translational and rotational part
of an element in~$\widetilde{\E}_2$ given in Section~\ref{section:deform}.
Elements with rotational part in $2\pi{\mathbb Z}$ are naturally
identified with translations of~${\mathbb R}^3$.

\begin{lem}
\label{lem:commute}
Let $s,t\in\pi$ be commuting elements of infinite order which generate
a (normal) subgroup ${\mathbb Z}\oplus{\mathbb Z}$ of~$\pi$. Then
$\overline{\rho}(s)$ and $\overline{\rho}(t)$ are translations
of~${\mathbb R}^2$, hence $\rho (s)$ and $\rho (t)$ are translations of
${\mathbb R}^3$.
\end{lem}

Normality of the subgroup generated by $s$ and $t$
will not be used in the following proof, but it will hold in all cases where we
apply this lemma.

\vspace{2mm}

\noindent {\em Proof.}
Write
\[ \rho (s)=\left(\left(\begin{array}{cc}x_s\\y_s\end{array}\right) ,
\theta_s\right).\]
We want to show that $\theta_s\in 2\pi {\mathbb Z}$. Arguing by contradiction,
we assume $\theta_s\not\in 2\pi {\mathbb Z}$, so $\overline{\rho}(s)$
is not a translation. Then $\overline{\rho}(s)$ is rotation about
a unique fixed point $z_s\in{\mathbb R}^2$. Since $s$ and $t$ commute,
we get
\[ \overline{\rho}(s)\overline{\rho}(t) (z_s)=\overline{\rho}(t)
\overline{\rho}(s) (z_s)=\overline{\rho}(t) (z_s),\]
hence $\overline{\rho}(t)(z_s)=z_s$. So $\overline{\rho}(s),
\overline{\rho}(t)$ are rotations about $z_s\in{\mathbb R}^2$ with angle
$\theta_s,\theta_t$, respectively.

If $\theta_t/\theta_s\in{\mathbb Q}$, then $m\theta_s+n\theta_t=0$
for suitable $m,n\in{\mathbb Z}\setminus\{ 0\}$. This implies that
$\rho (s^mt^n)$ is the identity element in~$\widetilde{\E}_2$, so
$\rho$ would not be faithful.

If $\theta_t/\theta_s\not\in{\mathbb Q}$, then, as is well known,
for any $\varepsilon >0$ one can find $m,n\in{\mathbb Z}$ such that
\[ 0<|m\theta_s+n\theta_t|<\varepsilon .\]
So in this case $\rho$ would not have a discrete image.

Since $\rho\in\R (\pi ,\widetilde{E}_2)$ is, by definition, faithful
and with discrete image, the assumption $\theta_s\not\in 2\pi{\mathbb Z}$
leads to a contradiction.\hfill $\Box$

\hspace{2mm}

To state the geometric description of the deformation spaces we recall
some standard notation. As in Section~\ref{section:sl} we write
\[ {\mathbb H}^2=\left\{ \tau\in{\mathbb C}\co\mbox{\rm Im}\,\tau >0\right\} \]
for the upper half-plane, on which $\SL_2{\mathbb Z}$ acts from
the left by fractional
linear transformations. Write ${\mathbb H}^2\times T^2$ for the trivial
$T^2$-bundle over ${\mathbb H}^2$ defined as a quotient
\[ ({\mathbb H}^2\times {\mathbb C})/\sim ,\]
where $(\tau ,z)\sim (\tau ',z')$ if and only if $\tau =\tau '$ and
$z-z'$ lies in the lattice $\langle 1,\tau \rangle\subset{\mathbb C}$
generated over ${\mathbb Z}$ by 1 and~$\tau$. This topologically
trivial $T^2$-bundle is the `universal elliptic curve': a bundle
over the Teichm\"uller space ${\mathbb H}^2$ of elliptic curves,
where the fibre over a point $\tau\in{\mathbb H}^2$ is exactly the
elliptic curve corresponding to that element of Teichm\"uller
space. Cf.\ Section~\ref{section:complex} for the analogous construction
of universal families for taut contact circles on $3$-manifolds.

The natural diagonal action of $\SL_2{\mathbb Z}$ on ${\mathbb H}^2\times
T^2$ is given by
\[ \SL_2{\mathbb Z}\ni\left(\begin{array}{cc}a&b\\c&d\end{array}\right)\co
(\tau ,z\;\mbox{\rm mod}\;\langle 1,\tau\rangle )\longmapsto
\left( \frac{a\tau +b}{c\tau +d},\frac{z}{c\tau +d}\; \mbox{\rm mod}\;
\langle 1, \frac{a\tau +b}{c\tau +d}\rangle \right) .\]

The following theorem describes the moduli spaces of taut contact circles
on the five left-quotients of~$\widetilde{\E}_2$. The subsequent proof
will also contain a description of the corresponding
Teichm\"uller spaces.

\begin{thm}
\label{thm:deform-e}
Let $M_k$ be the $T^2$-bundle over $S^1$ with periodic monodromy
of period~$k\in\{ 1,2,3,4,6\}$. The moduli spaces $\M (M_k)$ of taut
contact circles are given in the following table:

\begin{center}
\begin{tabular}{c|c}
$k$ & $\M (M_k)$\\ \hline
$1$ & ${\mathbb N}\times\SL_2{\mathbb Z}\backslash ({\mathbb H}^2\times
     T^2)$\rule{0cm}{.6cm}\\
$2$ & $\{ r\in{\mathbb N}\co r\equiv 1\;\mbox{\rm mod}\; 2\} \times
     \PSL_2{\mathbb Z}\backslash {\mathbb H}^2$\rule{0cm}{.6cm}\\
$3$ & $\{ r\in{\mathbb N}\co r\equiv \pm 1\;\mbox{\rm mod}\;
     3\}$\rule{0cm}{.6cm}\\
$4$ & $\{ r\in{\mathbb N}\co r\equiv \pm 1\;\mbox{\rm mod}\;
     4\}$\rule{0cm}{.6cm}\\
$6$ & $\{ r\in{\mathbb N}\co r\equiv \pm 1\;\mbox{\rm mod}\;
     6\}$\rule{0cm}{.6cm}\\
\end{tabular}
\end{center}
\end{thm}

\noindent {\em Proof.}
(i) For $M_1=T^3$, Lemma~\ref{lem:commute} applies to show that $\rho (u)$
is a translation for all $u\in\pi\cong{\mathbb Z}^3$, hence
$\R ({\mathbb Z}^3,\widetilde{\E}_2)\subset\R ({\mathbb Z}^3,
{\mathbb C}\times 2\pi {\mathbb Z})$ (after identifying the translational
parts of elements in $\widetilde{\E}_2$ with~${\mathbb C}$).
Observe that the
action $\rho\mapsto\rho^w$ of $w\in\widetilde{\E}_2$ by conjugation on
$\rho$ is as follows:
If $w$ is a translation
this action is trivial; if the translational
part of $w$ is~$0$, then $\rho^w(u)$ is obtained by rotating the translational
part of $\rho (u)$ by an angle equal (mod~$2\pi)$ to the rotational part
of~$w$. This implies that the left-action of $\Inn (\widetilde{\E}_2)$
and the right-action of ${\mathbb R}^+$ on $\R ({\mathbb Z}^3,{\mathbb C}
\times 2\pi{\mathbb Z})$ together constitute the standard
${\mathbb C}^*$-action on~${\mathbb C}^3$. Hence
\begin{eqnarray*}
\lefteqn{\T (T^3)=\left\{ ([z_1:z_2:z_3],(r_1,r_2,r_3))\in{\mathbb C}P^2
\times{\mathbb Z}^3\co\right. } \\
 & & \left. \mbox{\rm rank}_{\mathbb R} ((z_1,r_1),(z_2,r_2),
(z_3,r_3))=3\right\} .
\end{eqnarray*}
Notice that the points in $\T (T^3)$ with fixed coordinates
$(r_1,r_2,r_3)\in{\mathbb Z}^3$ form the complement of a real
hypersurface in~${\mathbb C}P^2$. For $(r_1,r_2,r_3)=(0,0,1)$,
the $z_j$-coordinates of a point in $\T (T^3)$ are of the form
$[1:z_2:z_3]$ with $z_2\not\in{\mathbb R}$. Hence, topologically
$\T (T^3)$ consists of a disjoint union of copies of~${\mathbb R}^4$.

We now turn to $\M (T^3)=\T (T^3)/\Out ({\mathbb Z}^3)$.
Given $\rho\in\R ({\mathbb Z}^3,\widetilde{E}_2)$, one can find a basis
$a_1,a_2,a_3$ for ${\mathbb Z}^3$ such that the rotational parts
of $\rho (a_1)$ and $\rho (a_2)$ are zero, and their translational parts
form a {\em negative} basis for~${\mathbb R}^2={\mathbb C}$ (with respect
to the complex orientation), and such that the rotational
part of $\rho (a_3)$ is equal to $2\pi r$ with $r\in {\mathbb N}$. Here
$r$ is the greatest common divisor of the rotational parts$/2\pi$ of
the $\rho (u)$, $u\in {\mathbb Z}^3$ (all of which lie in $2\pi{\mathbb Z}$).

The right action of $\Out ({\mathbb Z}^3)=\GL_3{\mathbb Z}$ allows us to
make arbitrary basis changes in ${\mathbb Z}^3$, so we can fix a basis
$e_1,e_2,e_3$ and assume that our representative $\rho$ of $[\rho ]\in
\M (T^3)$ maps this fixed basis as described above for $a_1,a_2,a_3$.
Any further action of $\Out ({\mathbb Z}^3)$ on special representatives
of this kind is then given by the subgroup ${\mathcal B}\subset
\SL_3{\mathbb Z}$ consisting of elements which stabilise
the subgroup ${\mathbb Z}^2\oplus 0\subset {\mathbb Z}^3$
generated by $e_1$ and~$e_2$. That is, elements of ${\mathcal B}$
are of the form
\[ \left(\begin{array}{ccc}
b_{11} & b_{12} & *\\
b_{21} & b_{22} & *\\
0 & 0 & 1\end{array}\right) ,\]
with $B=(b_{ij})\in\SL_2{\mathbb Z}$ (here $\SL_3{\mathbb Z}$ acts on
${\mathbb Z}^3$ from the left).

Furthermore, the ${\mathbb C}^*$-action on $\R ({\mathbb Z}^3,
\widetilde{\E}_2)$
allows us to assume that our representative $\rho$ is as follows:
\begin{eqnarray*}
\rho (e_1) & = & (\tau ,0),\;\;\tau\in{\mathbb H}^2,\\
\rho (e_2) & = & (1,0),\\
\rho (e_3) & = & (z,2\pi r),\;\; r\in {\mathbb N}.
\end{eqnarray*}
We call this a representation in {\em standard form}.
Thus, writing $\R ''({\mathbb Z}^3,\widetilde{\E}_2)$ for the space of
special representations as just described, we get
\[ \M (T^3)=\R ''({\mathbb Z}^3,\widetilde{\E}_2)/{\mathcal B},\]
and this translates into the description given in the theorem, in complete
analogy with the usual description of moduli space for surfaces
of genus~1. The only difference is due to the additional parameters
$r$ and $z$. The fact that $-\mbox{\rm id}\in\SL_2{\mathbb Z}$ acts
non-trivially on $z$ means that we cannot pass from $\SL_2{\mathbb Z}$
to~$\PSL_2{\mathbb Z}$. Notice that away from the fixed points of
the $\SL_2{\mathbb Z}$-action on~${\mathbb H}^2$, the space
$\SL_2{\mathbb Z}\backslash ({\mathbb H}^2\times T^2)$ is the universal
elliptic curve over the moduli space $\SL_2{\mathbb Z}\backslash
{\mathbb H}^2$ of elliptic curves; over the fixed points the fibre
is the quotient of an elliptic curve by a finite group action. 

\vspace{2mm}

(ii) As we pass to higher values of~$k$, the symmetries coming from
a non-trivial monodromy mean that we lose the freedom in the
parameters $z$ and~$\tau$. This accounts for the loss in
dimension of the corresponding moduli spaces. In all these remaining
cases $k\in\{ 2,3,4,6\}$ we fix a presentation of $\pi =\pi_{(k)}$
as at the beginning of this section. Let $\rho$ be an element of
$\R (\pi ,\widetilde{\E}_2)$. From Lemma~\ref{lem:commute} we know
that $\rho (s)$ and $\rho (t)$ are translations. The element $\rho (b)$,
on the other hand, cannot be a translation; otherwise it would commute with
$\rho (s)$ and $\rho (t)$, and $\rho$ would not be faithful.

Conjugation of $\rho$ by a translation in $\widetilde{\E}_2$ does not
affect translations, and it allows to change the fixed point in
${\mathbb R}^2={\mathbb C}$ of a rotation. So we can choose a representative
$\rho$ for $[\rho ]\in\T (M_k)$, $k\in\{ 2,3,4,6\}$, in such a way
that $\overline{\rho}(b)$ is a rotation about $0\in{\mathbb C}$.

Since $b^k$ does commute with $s$ and~$t$, Lemma~\ref{lem:commute}
shows that $\overline{\rho}(b^k)$ is a translation. But it fixes $0\in
{\mathbb C}$, so it must be trivial. Hence $\rho (b)=(0,2\pi r/k)$
with $r\in{\mathbb Z}-\{ 0\}$.

As earlier we write $\theta_u$ for the rotational part of
$\rho (u)\in\widetilde{\E}_2$, $u\in\pi$. We have $\theta_{uv}=
\theta_u+\theta_v$, so the relations in $\pi$ yield the following
equations:
\begin{eqnarray*}
\theta_s & = & \alpha\theta_s+\gamma\theta_t,\\
\theta_t & = & \beta\theta_s+\delta\theta_t,
\end{eqnarray*}
or as a matrix equation,
\[ \left(\begin{array}{cc}\alpha -1 &\gamma\\ \beta & \delta -1
\end{array}\right) \left( \begin{array}{c}\theta_s\\ \theta_t
\end{array} \right ) = \left( \begin{array}{c}0\\0
\end{array} \right ) .\]
The determinant of this $(2\times 2)$-matrix is
\[ (\alpha -1)(\delta -1)-\beta\gamma = \alpha\delta -\beta\gamma -
(\alpha +\delta )+1 = 2-\mbox{\rm trace}\, A.\]
For a periodic matrix $A\in\SL_2{\mathbb Z}$ different from the identity
matrix we have $\mbox{\rm trace}\, A\neq 2$, so the only solution to the
equation above is $\theta_s=\theta_t=0$. Therefore, any element of
$\T (M_k)$ has a representative of the form
\begin{eqnarray*}
\rho (s) & = & (1,0),\\
\rho (t) & = & (\tau ,0),\;\;\tau\in{\mathbb C}-{\mathbb R},\\
\rho (b) & = & (0,2\pi r/k),\;\; r\in{\mathbb Z}-\{ 0\}.
\end{eqnarray*}
(As in (i), rotation of the translational parts can be effected by
conjugation with a rotation $w\in\widetilde{\E}_2$.)

\vspace{2mm}

(iii) Now consider the different values of $k$ in turn, beginning
with~$k=2$: To get a representation for $\pi_{(2)}$ it is necessary and
sufficient to have $r$ odd, so
\[ \T (M_2)=\{ r\in {\mathbb Z}\co r\equiv 1\;\mbox{\rm mod}\; 2\}
   \times ({\mathbb C}-{\mathbb R}).\]

The group $\pi_{(2)}$ admits automorphisms which send any of the generators
$s,t,b$ to their inverse and fixes the others. So in $\M (M_2)=
\T (M_2)/\Out (\pi_{(2)})$ we may impose on the representative $\rho$
for a given class the further conditions $\mbox{\rm Im}\,\tau >0$ and
$r\in{\mathbb N}$.

Any further action of $\Out (\pi_{(2)})$ on representations of this kind must
(and can) be of the form
\[ s\longmapsto s^a t^b,\;\;\; t\longmapsto s^c t^d\]
for arbitrary $\left( \begin{array}{cc}a&b\\c&d\end{array}\right)\in\SL_2
{\mathbb Z}$, where the action of $-\mbox{\rm id}\in\SL_2{\mathbb Z}$ is
trivial modulo $\Inn (\widetilde{\E}_2)$. It follows
that $\M (M_2)$ is as claimed.

\vspace{2mm}

(iv) $k=3,4,6$: The arguments in these three cases are completely analogous
to one another, and we only present the one for~$k=3$. The presentation
of $\pi_{(3)}$ is
\[ \pi_{(3)}=\left\{ s,t,b\co st=ts,\; bsb^{-1}=t,\; btb^{-1}=s^{-1}t^{-1}
\right\} .\]
Since the monodromy has order~$3$, we must have $r\equiv\pm 1$ mod~$3$
in the description of~$\rho$. Then the relation $bsb^{-1}=t$ forces
$\tau =\exp (2\pi i\, r/3)$. This is consistent with $btb^{-1}=s^{-1}
t^{-1}$, indeed, the given $\tau$ satisfies $\tau^2=-1-\tau$. Hence
\[ \T (M_3)=\{ r\in{\mathbb Z}\co r\equiv\pm 1\;\mbox{\rm mod}\; 3\} .\]
The group $\pi_{(3)}$ admits an automorphism defined by
\[ s\longmapsto s,\;\; t\longmapsto s^{-1}t^{-1},\;\; b\longmapsto b^{-1},\]
as can also be seen from the geometry of the torus bundle, and there
are no other non-trivial outer automorphisms to consider. Thus
\[ \M (M_3)=\T (M_3)/\pm =\{ r\in {\mathbb N}\co r\equiv\pm 1\;
\mbox{\rm mod}\; 3\} .\]

This completes the proof of Theorem~\ref{thm:deform-e}.\hfill $\Box$
\section{$\SU (2)$-geometry}
\label{section:s}
In this section we return to the study of taut contact circles on
left-quotients $M$ of~$\SU (2)$. Our aim is to obtain information about
the Teichm\"uller spaces $\T (M)$ of taut contact circles on these
manifolds. Since $\M (M)$ is known (Theo\-rem~\ref{thm:deform-s}) and
\[ \M (M)=\T (M)/(\Diff (M)/\Diff_0(M)), \]
we need to investigate the action of the {\bf homeotopy group}
\[ \Lambda (M)=\Diff (M)/\Diff_0(M) \]
on taut contact circles. Notice that, as a set, this homeotopy group
can be identified with the set $\pi_0(\Diff (M))$ of path-connected
components of $\Diff (M)$.

First we recall the list of finite subgroups of~$\SU (2)$,
see~\cite[Thm.~2.6.7]{wolf84}. We also fix a particular representation
$\rho$ of these groups in~$\SU (2)$. Thus we do not distinguish between the
abstract group $\pi$ and its image $\Gamma =\rho (\pi )\subset\SU (2)$,
and we identify $u\in\pi$ both with $\rho (u)\in\Gamma$ and
its action $L_{\rho (u)}$ by left multiplication. For our purposes this
is justified since, by the cited theorem, isomorphic finite subgroups
of $\SU (2)$ are actually conjugate in~$\SU (2)$.

We identify $\SU (2)$ with the unit quaternions $S^3\subset{\mathbb H}$
by the group isomorphism
\[ \left( \begin{array}{rr}
a_0+ia_1 & b_0+ib_1\\ -b_0+ib_1 & a_0-ia_1
\end{array} \right) \longmapsto a_0+ia_1+jb_0+kb_1. \]

\begin{itemize}
\item The {\bf cyclic group} $C_m$ of order~$m$, generated by $x=\cos (2\pi /m)
+i\sin (2\pi /m)$.
\item The {\bf binary dihedral group} of order~$4n$,
\[ D_{4n}^*=\left\{ x,y\co \, x^2=(xy)^2= y^n\right\} ,\]
generated by $x=i$ and $y=\cos (\pi /n)+j\sin (\pi /n)$. The group
\[ D_8^*=\left\{ \pm 1,\pm i,\pm j,\pm k\right\} \]
is also called the {\bf quaternion group} $Q_8$.
\item The {\bf binary tetrahedral group} of order~$24$,
\[ T^*=Q_8^{x,y}\rtimes C_3^z,\]
generated by $x=i$, $y=j$, and $z=-(1+i+j+k)/2$. Here $C_3$ acts on
$Q_8$ by $zxz^{-1}=y$ and $zyz^{-1}=xy$.
\item The {\bf binary octahedral group} of order~$48$,
\[ O^*=T^*\rtimes C_4^w/(w^2=x^2), \]
with $w=(i-k)/\sqrt{2}$. The action of $C_4$ on $T^*$ is given by
$wxw^{-1}=yx$, $wyw^{-1}=y^{-1}$, and $wzw^{-1}=z^{-1}$.
\item the {\bf binary icosahedral group} of order $120$ with presentation
\[ I^*=\left\{ a,b\co \, a^2=(ab)^3=b^5, a^4=1\right\} .\]
\end{itemize}

\begin{rem}
{\rm The groups $T^*$ and $O^*$ also admit the more concise presentations
\[ \bigl\{ a,b\co \, a^2=(ab)^3=b^{\beta}, a^4=1\bigr\} ,\; \beta =3,4,\]
but the presentations above are more useful for our purposes. For $T^*$,
for instance, the isomorphism between the two presentations is given
by $a\mapsto x$ and $b\mapsto yz$.}
\end{rem}

Let $\oph$ be a diffeomorphism of $M=\Gamma\backslash\SU (2)$. This lifts
to a diffeomorphism $\varphi$ of~$\SU (2)$, which defines an
automorphism $\vartheta_{\oph}$ of $\Gamma$ by
\[ \varphi\circ u\circ\varphi^{-1}=\vartheta_{\oph}(u)
\;\; \forall u\in\Gamma .\]
Different lifts $\varphi$ of $\oph$ differ by a deck transformation
(i.e.\ an element of~$\Gamma$), so the equivalence class
$[\vartheta_{\oph}]\in\Out (\Gamma )$ is determined by~$\oph$, and it is easy
to check that $\oph\mapsto \vartheta_{\oph}$ defines a homomorphism.

If $\oph_t$ is an isotopy from $\oph_0=\oph$ to $\oph_1=\mbox{\rm id}$,
there is a lifted isotopy $\varphi_t$ from $\varphi_0=\varphi$ to
a deck transformation $\varphi_1=w\in\Gamma$. During this isotopy the
automorphism $\vartheta_{\oph}$ of the discrete group $\Gamma$ is forced to
stay fixed. So there is a well-defined homomorphism
\[ \begin{array}{rccc}
\chi\co & \Lambda (M) & \longrightarrow & \Out (\Gamma )\\
  & [\oph ] & \longmapsto & [\vartheta_{\oph}].
\end{array} \]

We first consider the case that $M=\Gamma\backslash\SU (2)$ is a quotient
of $\SU (2)$ under one of the non-abelian subgroups $\Gamma$ of~$\SU (2)$.
We know that in this case the moduli space $\M (M)$ consists of a single
point. Since $\M (M)=\T (M)/\Lambda (M)$, it follows that $\T (M)$ is
a discrete set on which $\Lambda (M)$ acts transitively. So the order
$| \T (M)|$ of this set is equal to $| \Lambda (M)/\Lambda_0(M)|$, where
$\Lambda_0(M)$ denotes the subgroup of $\Lambda (M)$ stabilising
a particular element of~$\T (M)$, say the one represented by
$\omega_1+i\omega_2=z_1dz_2-z_2dz_1$. In the sequel we are going
to determine~$| \Lambda_0(M)|$, and with that information we obtain the
following theorem.

\begin{thm}
\label{thm:teich-s1}
Let $M=\Gamma\backslash\SU (2)$ be a left-quotient of $\SU (2)$ under a
discrete, cocompact, non-abelian subgroup~$\Gamma$. Then $\T (M)$ is a
discrete set, and the order of this set is given in the following table,
where we write $|\Lambda |=|\Lambda (M)|$ for short:
\begin{center}
\begin{tabular}{c||c|c|c|c|c|c}
$\Gamma$ & $Q_8$ & $D_{8n+4}^*$ & $D_{8n+8}^*$ &
$T^*$ & $O^*$ & $I^*$\\ \hline
$|\T (M)|$ & $|\Lambda |/6$ & $|\Lambda |/2$ & $|\Lambda |$ &
$|\Lambda |/2$ & $|\Lambda |$ & $|\Lambda |$  \rule{0cm}{.6cm}\\
\end{tabular}
\end{center}
\end{thm}

\begin{rem}
{\rm We are not aware of any information about $|\Lambda (M)|$ for the
non-abelian groups~$\Gamma$.}
\end{rem}

\begin{lem}
\label{lem:nonab}
Let $M$ be a non-abelian left-quotient of~$\SU (2)$. Then
\[ \chi |_{\Lambda_0(M)}\co \Lambda_0(M)\longrightarrow \Out (\Gamma ) \]
is a monomorphism onto the subgroup $\Out_0(\Gamma )$ of outer
automorphisms of $\Gamma\subset\SU (2)$ that are given by conjugation
with a matrix in~$\SL_2{\mathbb C}$.
\end{lem}

\noindent {\em Proof.}
Let $[\oph ]$ be an element of $\Lambda_0(M)$. This means that
$\oph$ preserves the homothety class of $\omega_1+i\omega_2=
z_1dz_2-z_2dz_1$ up to a diffeomorphism isotopic to the identity.
Hence, without changing the class $[\oph ]$ we may actually assume that
$\oph$ fixes this homothety class. Moreover, the rotation of
$(\omega_1,\omega_2)$ can be effected (on~$S^3$) by the flow of
the left-invariant vector field $X$ defined by $\omega_1(X)=\omega_2(X)=0$
and $d\omega_1(X,.)=\omega_1$ (this is simply the Hopf vector field).
This flow descends to~$M$, so we may assume that $\oph$ fixes
the conformal class of $(\omega_1,\omega_2)$, say $\oph^*(\omega_1,
\omega_2)=\ov (\omega_1,\omega_2)$ with $\ov\co M\rightarrow
{\mathbb R}^+$. Lift $\oph$ to a diffeomorphism $\varphi$ of $S^3$
and $\ov$ to a function $v$ on~$S^3$. From the construction
in \cite[Section~3]{gego95} it follows that the diffeomorphism
\[ \begin{array}{rccc}
\phi\co & S^3\times{\mathbb R} & \longrightarrow & S^3\times{\mathbb R}\\
  & (x,t) & \longmapsto & (\varphi (x),t-\log v(x))
\end{array} \]
is an ${\mathbb R}$- and $\Gamma$-equivariant holomorphic automorphism
of ${\mathbb C}^2-\{ (0,0)\}$ under the natural identification
of that space with $S^3\times{\mathbb R}$, and it preserves the
holomorphic $1$-form $\omega_1+i\omega_2$. This $\phi$ extends to
an automorphism of ${\mathbb C}^2$ fixing the origin. The
${\mathbb R}$-equivariance implies that the holomorphic differential
of $\phi$ is bounded, and hence constant. Therefore $\phi$ has to be
a linear automorphism, and the fact that it preserves $d(\omega_1+i
\omega_2)=dz_1\wedge dz_2$ forces $\phi\in\SL_2{\mathbb C}$.
The $\Gamma$-equivariance of $\phi$ takes the form
\[ \phi\circ u\circ\phi^{-1}=\vartheta_{\oph}(u)\;\;\forall u\in\Gamma ,\]
where the left-hand side may now be read as matrix multiplication
in~$\SL_2{\mathbb C}$.

Conversely, if we start with an element of $\SL_2{\mathbb C}$ that
normalises~$\Gamma$, then it induces an ${\mathbb R}$-equivariant
diffeomorphism of $\Gamma\backslash ({\mathbb C}^2-\{ (0,0)\})
\cong M\times {\mathbb R}$ that preserves $\omega_1+i\omega_2=
z_1dz_2-z_2dz_1$, and hence a diffeomorphism of $M$ preserving the conformal
class of $(\omega_1,\omega_2)$. This proves that $\chi (\Lambda_0(M))
=\Out_0(\Gamma )$.

It remains to prove injectivity of $\chi |_{\Lambda_0(M)}$. Since the
image $\chi (\Lambda_0(M))$ consists of automorphisms given by
conjugation with an element $\phi\in\SL_2{\mathbb C}$, it suffices to
show that the condition
\[ \phi\circ u\circ\phi^{-1}=wuw^{-1}\;\;\forall u\in\Gamma , \]
with $\phi\in\SL_2{\mathbb C}$ and some $w\in\Gamma$, forces $\phi
\in\Gamma$. But this follows easily from the fact that the
${\mathbb C}$-linear span of any non-abelian subgroup $\Gamma\subset\SU (2)$
consists of all complex $(2\times 2)$-matrices, so the preceding condition
(which is ${\mathbb C}$-linear in~$u$) does in fact imply $\phi =w$.
This concludes the proof of the lemma. \hfill$\Box$

\vspace{2mm}

Together with the following proposition, this lemma constitutes
a proof of Theorem~\ref{thm:teich-s1}.

\begin{prop}
\label{prop:out-nonab}
The groups $\Out_0(\Gamma )$ are as given in the following table:
\begin{center}
\begin{tabular}{c||c|c|c|c|c|c}
$\Gamma$ & $Q_8$ & $D_{8n+4}^*$ & $D_{8n+8}^*$ &
$T^*$ & $O^*$ & $I^*$\\ \hline
$\Out_0 (\Gamma )$ & $S_3$ & $C_2$ & $1$ & $C_2$ & $1$ & $1$
\rule{0cm}{.6cm}\\
\end{tabular}
\end{center}
\end{prop}

\noindent {\em Proof.}
(i) $Q_8$: The group $\Aut (Q_8)$ contains as a normal subgroup the
Klein $4$-group~$V_4$, corresponding to the sign changes of $i$
and~$j$ (the sign of $k$ is then determined). These are inner automorphisms
of~$Q_8$, given by conjugation with $i$, $j$, or~$k$. The quotient
group $\Aut (Q_8)/V_4$ is isomorphic to the symmetric group~$S_3$,
acting by permutations of $i$, $j$, and~$k$. The transposition $(ij)$
is given by conjugation with the unit quaternion $(1+k)/\sqrt{2}$,
\[ \frac{1+k}{\sqrt{2}}\, i\, \frac{1-k}{\sqrt{2}}=j.\]
The $3$-cycle $(ijk)=(jk)(ik)$ is given by conjugation with
\[ \frac{1+i}{\sqrt{2}}\frac{1-j}{\sqrt{2}}=\frac{1+i-j-k}{2}. \]
These are not inner automorphisms, but they are generated by conjugation
with unit quaternions, and hence elements of $\SU (2)\subset
\SL_2{\mathbb C}$ under our identification of the unit quaternions
with~$\SU (2)$. It follows that
\[ \Out (Q_8)=\Out_0(Q_8)=S_3.\]

\vspace{2mm}

(ii) $D_{4n}^*$, $n\geq 3$: From the relation $x^2=(xy)^2=y^n$ one
obtains $x^{-1}yx=y^{-1}$, hence
\[ x^2=x^{-1}x^2x=x^{-1}y^nx=y^{-n}=x^{-2},\]
so $x$ is of order~$4$ and $y$ of order~$2n$. Since $yx=xy^{-1}$,
a complete list of elements of $D_{4n}^*$ is given by
\[ \left\{ 1,y,\ldots ,y^{2n-1},x,xy,\ldots ,xy^{2n-1}\right\} .\]
The order of each element of the form $xy^k$ is equal to~$4$. It follows
that automorphisms of $D_{4n}^*$, $n\geq 3$, can only be of the form
$x\mapsto xy^l$, $y\mapsto y^k$ with $\gcd (k,2n)=1$. The assumption
that such an automorphism is given by conjugating $D_{4n}^*\subset
\SU (2)$ with an element $\phi\in\SL_2{\mathbb C}$ leads, via some
simple algebra, to the following list of possibilities:
\[ \begin{array}{ll}
x\mapsto x, & y\mapsto y^{2n-1};\\
x\mapsto xy^n=x^{-1}, & y\mapsto y;\\
x\mapsto xy^n=x^{-1}, & y\mapsto y^{2n-1}.
\end{array} \]
If we regard $D_{4n}^*$ as a subgroup of the unit quaternions, these
automorphisms are induced by conjugation with $i$, $j$ and~$k$,
respectively. So the first one is always inner; the second and third one
are inner if and only if $n$ is even (and they always differ by
an inner automorphism).

\vspace{2mm}

(iii) $T^*$: This group contains $Q_8$ as a normal subgroup, and it
is a straightforward check from the explicit presentations that
the non-trivial elements of $Q_8\subset T^*$ are characterised as
exactly those elements of $T^*$ that have order~$4$. It follows that
any automorphism of $T^*$ induces an automorphism of~$Q_8$. So we can
refer to the results in~(i).

The automorphism of $Q_8$ given by conjugation with $(1+k)/\sqrt{2}$
also induces one of~$T^*$, and it is not inner. Conjugation with
$(1+i-j-k)/2=zy$, on the other hand, is an inner automorphism.
This implies $\Out_0(T^*)=S_3/A_3=C_2$.

\vspace{2mm}

(iv) $O^*$: This group contains a unique copy of~$T^*$, so any
automorphism of $O^*$ will induce one of~$T^*$. Since
$xwz=(1+k)/\sqrt{2}$, the only possible outer automorphism of
$T^*$ given by conjugation with an element of $\SL_2{\mathbb C}$
actually comes from an inner automorphism of~$O^*$.
We infer that $\Out_0(O^*)$ is trivial.

\vspace{2mm}

(v) $I^*$: It is known that $\Out (I^*)=C_2$,
cf.~\cite[p.~195]{wolf84}. Moreover, $I^*$ has exactly two
non-equivalent irreducible complex fixed-point free representations,
both of degree~$2$, which differ by the non-trivial outer
automorphism, see~\cite[Lemma~7.1.7]{wolf84}.
So this outer automorphism cannot come from
conjugation with an element in $\SL_2{\mathbb C}$.

This concludes the proof of the proposition.\hfill $\Box$

\vspace{2mm}

We now turn to the lens spaces $L(m,m-1)$. Use the notation
$\M_1,\M_2$ as in Theorem~\ref{thm:deform-s} and write
$\widetilde{\M}_1$ for the complex slab $\{ a\in{\mathbb C}\co
0<\mbox{\rm Re}(a)<1\}$.

\begin{thm}
\label{thm:teich-s2}
For $M=L(m,m-1)$, the Teichm\"uller spaces $\T (M)$ are as follows:
\begin{center}
\begin{tabular}{c|c}
$m$ & $\T (L(m,m-1))$\\ \hline
$1,2$ & $\M_1\sqcup \M_1\sqcup\M_2\sqcup\M_2$ \rule{0cm}{.6cm} \\
$\geq 3$ & $\widetilde{\M}_1\sqcup\M_2\sqcup\M_2$ \rule{0cm}{.6cm}\\
\end{tabular}
\end{center}
\end{thm}

\noindent {\em Proof.}
(i) If $m=1$ or $2$ (i.e.\ $M=S^3$ or ${\mathbb R}P^3$), the
identification of elements in $\C(L(m,m-1))$ corresponding to
the parameter values $a,1-a$ (in the continuous family
of taut contact circles described in Section~\ref{section:deform})
can be effected by the diffeomorphism of $L(m,m-1)$ induced by
$(z_1,z_2)\mapsto (z_2,-z_1)$. This diffeomorphism lies in
$\Diff_0(M)$ in both cases. Furthermore, we have
$\Lambda (M)=C_2$, generated by the orientation reversing
diffeomorphism $(z_1,z_2)\mapsto (\overline{z}_1,z_2)$. So
$\Lambda (M)$ acts non-trivially on $\T (M)$, and we conclude
that $\T (M)$ consists of two disjoint copies of~$\M (M)$.

\vspace{2mm}

(ii) If $m\geq 3$, we still have $\Lambda (M)=C_2$,
see~\cite{bona83},~\cite{horu85}.
But the diffeomorphism $(z_1,z_2)\mapsto (z_2,-z_1)$ conjugates the
action of $x=\cos (2\pi /m)+i\sin (2\pi /m)$ to that of~$\overline{x}$.
Since $x$ and $\overline{x}$ are not conjugate in $C_m$ for~$m\geq 3$,
the induced diffeomorphism of $L(m,m-1)$ cannot be isotopic to the
identity. (In particular, this orientation preserving diffeomorphism
defines the non-trivial element in~$\Lambda (M)$, and $M$ does not
admit any orientation reversing diffeomorphism.)
As in the proof of Lemma~\ref{lem:nonab} we see that
a diffeomorphism of $L(m,m-1)$ that preserves the homothety class
of one of the standard models for a taut contact circle, or sends
the model with parameter value~$a$ to that with parameter~$1-a$,
has to be induced by conjugation with an element $\phi\in\SL_2
{\mathbb C}$. This forces $\phi$ to be diagonal or anti-diagonal.
The latter happens exactly when we exchange $a$ and $1-a$ and the
automorphism of $C_m$ induced by $\phi$ is $x\mapsto\overline{x}$.
It follows that $\T (M)$ is as claimed. \hfill $\Box$
\section{The complex structure of Teichm\"uller space}
\label{section:complex}
For manifolds $M$ modelled on $\widetilde{\SL}_2$ we have seen in
Theorem~\ref{thm:deform-sl} that $\T (M)\rightarrow \T (O_M)$ is a trivial
principal ${\mathbb Z}^{2g}$-bundle, so we can equip $\T (M)$ with a
complex structure such that the projection onto $\T (O_M)$ becomes a local
biholomorphism, where the complex structure on $\T (O_M)$ is the one from
classical Teichm\"uller theory. For manifolds modelled on $\SU (2)$ or
$\widetilde{\E}_2$ we have described the moduli space $\M (M)$
explicitly as a complex space. Our aim in the present section is to show
that these complex structures on $\T (M)$ or $\M (M)$ are the `natural'
ones to consider.

For the geometry $\SU (2)$ this may be justified from the explicit description
of the contact circles corresponding to points in $\M_1$
(see Theorem~\ref{thm:deform-s}) in terms of the complex parameter~$a$. This
description hinges on our theory developed in~\cite{gego95} which relates
taut contact circles on~$M$ to complex structures on $M\times {\mathbb R}$.
Here, by contrast and more intrinsically, we show that one can associate
with any taut contact circle a complex analogue of the Godbillon-Vey
invariant, which essentially recovers the moduli parameter~$a$ in the case
of left-quotients of~$\SU (2)$.
That complex Godbillon-Vey invariant is in fact an
invariant of transversely conformal flows on $3$-manifolds; this aspect will be
pursued further in a forthcoming paper~\cite{gegoc}.

For the geometries $\widetilde{\SL}_2$ and $\widetilde{\E}_2$ we describe
universal families analogous to those of Bers~\cite{bers73} in classical
Teichm\"uller theory. Such a universal family consists of a holomorphic
fibration over $\T (M)$, where the fibre over a point $\sigma\in
\T (M)$ is the complex surface (diffeomorphic to $M\times S^1$) determined --
in the sense of~\cite{gego95} -- by the taut contact circle on $M$
corresponding to~$\sigma$. Ana\-logous universal families for Seifert
$4$-manifolds are described in~\cite{ue92}.

\vspace{2mm}

{\bf (1)} In the case of the geometry $\widetilde{\E}_2$ we only need to
consider $k=1$ (i.e.\ taut contact circles on the $3$-torus). For $k=2$ the
construction will be similar. For the other values of $k$ the
Teichm\"uller spaces are discrete sets, so geometrically the construction
of a universal family is not an issue, but it might be of some arithmetic
interest.

For the $3$-torus a universal family can best be defined over
${\mathcal R}''({\mathbb Z}^3,
\widetilde{\E}_2)$ (and it can then be pulled back to
$\T (T^3)$). Given $[\rho ]\in{\mathcal R}''({\mathbb Z}^3,
\widetilde{\E}_2)$, choose the unique representative $\rho$ in the standard
form described in Section~\ref{section:e}.
Then define a representation $\rho^{\mathbb C}\co
{\mathbb Z}^4\rightarrow {\mathbb C}^2$ as follows:
\begin{eqnarray*}
\rho^{\mathbb C}(e_1) & = & (\tau ,0),\\
\rho^{\mathbb C}(e_2) & = & (1,0),\\
\rho^{\mathbb C}(e_3) & = & (z,2\pi ir),\\
\rho^{\mathbb C}(e_4) & = & (0,1).
\end{eqnarray*}
Recall that the identification of $\widetilde{\E}_2\times \E^1$ (where
$\E^1$ denotes the Euclidean line) with
${\mathbb C}^2$ is given by $(z,\theta ,\lambda )\equiv (z,\lambda +
i\theta )=(z,w)$,
cf.~\cite[p.~192]{gego95}.

Now let ${\mathcal E}_{\mathbb C}$ be the quotient space of
${\mathcal R}''({\mathbb Z}^3,
\widetilde{\E}_2)\times {\mathbb C}^2$ under the equivalence relation
\begin{eqnarray*}
\lefteqn{([\rho_1],(z_1,w_1))\sim ([\rho_2],(z_2,w_2))
:\Longleftrightarrow}\\
  &  & [\rho_1]=[\rho_2]\;\mbox{\rm and}\;
       (z_2-z_1,w_2-w_1)\in\rho^{\mathbb C}_1({\mathbb Z}^4).
\end{eqnarray*}
We have a natural complex structure on
\[ {\mathcal R}''({\mathbb Z}^3, \widetilde{\E}_2)={\mathbb N}\times
{\mathbb H}^2\times {\mathbb C}.\]
Using this and
the identification of $[\rho ]$ with $(r,\tau ,z)$, the space
${\mathcal E}_{\mathbb C}$ inherits a complex structure such that the projection
\[\begin{array}{rccc} \Psi_{\mathbb C}\co &  {\mathcal E}_{\mathbb C} &
\longrightarrow &
{\mathcal R}''({\mathbb Z}^3, \widetilde{\E}_2)\\
  & \bigl[ [\rho ],(z,w)\bigr] & \longmapsto & [\rho ]
\end{array} \]
is holomorphic. The fibres of $\Psi_{\mathbb C}$ are complex tori,
and it is immediate
from the construction in~\cite{gego95} that $\Psi_{\mathbb C}^{-1}([\rho ])$
is biholomorphically equivalent to the complex torus naturally associated to
the taut contact circle on $T^3$ defined by~$\rho$ (with $t_0=1$ in the
notation of~\cite[p.~191]{gego95}). These considerations show that the
complex structure on $\R '' ({\mathbb Z}^3,\widetilde{\E}_2)$ is the
one naturally adapted to the study of taut contact circles within
the complex geometric framework developed in~\cite{gego95}.

It is therefore appropriate to define a universal family
${\mathcal E}={\mathcal E}(T^3)$ of taut
contact circles on $T^3$ as the quotient space of $\R ''({\mathbb Z}^3,
\widetilde{\E}_2)\times\widetilde{\E}_2$ under the analogue of the equivalence
relation above, with ${\mathbb C}^2=\widetilde{\E}_2\times
\E^1$ replaced by $\widetilde{\E}_2$ and
$\rho^{\mathbb C}$ by~$\rho$. Then
\[\begin{array}{rccc} \Psi\co &  {\mathcal E} & \longrightarrow &
{\mathcal R}''({\mathbb Z}^3, \widetilde{\E}_2)\\
  & \bigl[ [\rho ],(z,i\theta )\bigr] & \longmapsto & [\rho ]
\end{array} \]
is a smooth map such that $\Psi^{-1}([\rho ])$ is a copy of $T^3$ equipped
with the taut contact circle corresponding to~$\rho$. The identification of
$\Psi^{-1}([\rho ])$ with $T^3$ is determined, up to an element of
$\Diff_0(T^3)$, by the requirement that, for a fixed set of generators
of $\pi_1(T^3)$, the developing map of Section~\ref{section:Cartan} give the
representation $\rho$ in standard form.

The action of $\B$ on $\R ''({\mathbb Z}^3,\widetilde{\E}_2)$ can be extended
to an action on~${\mathcal E}$ as follows. Let $\rho$ be
a representative in standard form of an element $[\rho ]$ in
$\R ''({\mathbb Z}^3,\widetilde{\E}_2)$, and let $\vartheta\in\B$.
Then $\vartheta$ acts on
such a standard representative via
\[ \rho\longmapsto \mu(\rho\circ \vartheta ),\]
where $\mu =\mu (\rho ,\vartheta )\in {\mathbb C}^*$ is determined by the
condition that $\mu (\rho\circ \vartheta )$ be again in standard form.

The fixed points of this action are exactly those $[\rho ]$ for which
the $2$-torus $T^2={\mathbb C}^2/({\mathbb Z}\oplus {\mathbb Z}\tau)$ has
a symmetry, given by multiplication by some $\mu\in{\mathbb C}^*$, for
which $z\in T^2$ is a fixed point. So we obtain a fibre space
\[ {\mathcal E}/\B\longrightarrow \M (T^3) ,\]
where the fibres are copies of $M$ with a taut contact circle determined
by $[\rho ]\in\M (T^3)$, quotiented out by symmetries of the kind described.
The possible symmetries are given in the following table:

\vspace{2mm}

\begin{center}
\begin{tabular}{c|c|c}
$\mu$ & $\tau$ & $z$\\ \hline
$-1$ ($2$-fold) & arbitrary & $0,1/2,
\tau/2,1/2+\tau /2$\rule{0cm}{.7cm}\\
$\exp (2\pi i/3)$ ($3$-fold) & $\exp (2\pi i/6),\exp (2\pi i/3)$ &
    $1/2+i/(2\sqrt{3}), 1+i/\sqrt{3}$\rule{0cm}{.7cm}\\
$i$ ($4$-fold) & $i$ & $0, 1/2+i/2$\rule{0cm}{.7cm}\\
\end{tabular}
\end{center}

\vspace{3mm}

{\bf (2)} In the case of the geometry $\widetilde{\SL}_2$ we build a universal
family from the universal family of Bers~\cite{bers73} in analogy with
the description for Seifert $4$-manifolds given by Ue~\cite{ue92}.

According to \cite{bers73} there is a fibre space
\[ \F (O_M)=\left\{ (\osi ,z) \in \T (O_M)\times {\mathbb C}\co z\in D(\osi )
\right\} \]
over $\T (O_M)$, where $D(\osi )$ is a domain in ${\mathbb C}$ on which
$\pi^{\orb}=\pi_1^{\orb}(O_M)$ acts according to~$\osi$. With the complex
structure on $\T (O_M)$ from classical Teich\-m\"uller theory,
$D(\osi )$ and the group $\osi (\pi^{\orb})$ depend holomorphically
on $\osi$ in a sense explained in~\cite{bers73}. Holomorphic dependence
of $\osi (\pi^{\orb})$ on $\osi$ means that the map
$(\osi ,z)\mapsto (\osi ,\osi (\ou )(z))$ is holomorphic in
$\osi\in\T (O_M)$ and $z\in D(\osi )$
for each fixed $\ou\in\pi^{\orb}$.

So $\osi$ is interpreted as a representation of $\pi^{\orb}$ in the
holomorphic automorphisms of~$D(\osi )$,
and we have a biholomorphism $h_{\osi}\co {\mathbb H}^2\rightarrow
D(\osi )$ and a representation $\orh\in\R (\pi^{\orb},\PSL_2{\mathbb R})$
with $[\orh ]=\osi\in\T (O_M)$
such that there is a commutative diagram
\[
\begin{CD}
{\mathbb H}^2 @>\orh (\ou )>> {\mathbb H}^2\\
@Vh_{\osi}VV             @VVh_{\osi}V\\
D(\osi ) @>\osi (\ou )>> D(\osi ).
\end{CD}
\]
No statement is made about the dependence of $h_{\osi}$ on~$\osi$.
Taking the quotient of $D(\osi )$
by $\osi (\pi^{\orb})$ in the fibre of $\F (O_M)$ over $\osi$ gives
the holomorphic family ${\mathcal E}(O_M)\rightarrow \T(O_M)$ of complex
structures on~$O_M$.

To get the corresponding universal family of taut contact circles,
the idea is simply to construct the logarithmic differential
of the commutative diagram above.
To this end, identify $\widetilde{\SL}_2\times\E^1$ with ${\mathbb H}^2\times
{\mathbb C}$ with coordinates $(z,w)$,
so that $w=\lambda +i\theta$ corresponds to
$\log dz$, with $\lambda$ denoting the $\E^1$-coordinate.
Given ${\mathbf w}\in {\mathbb Z}^{2g}$, write $\sigma =\sigma_{\mathbf w}
=[\rho_{\mathbf w}]\in\T (M)$ for the corresponding lift of
$\osi =[\orh ]\in\T (O_M)$. Regarding ${\mathbb H}^2\times
{\mathbb C}$ and $D(\osi )\times{\mathbb C}$ as holomorphic tangent
bundles, it makes sense to define a biholomorphism
${\mathbb H}^2\times {\mathbb C}\rightarrow D(\osi )\times {\mathbb C}$ by
\[ h_{\sigma}(z,w) =\bigl(h_{\osi}(z),w+\log\Bigl(
\frac{\partial h_{\osi}}{\partial z}\Bigr)
(z) \bigr) .\]
Here the branch of the logarithm may be chosen arbitrarily (but
holomorphically in~$z$; this is possible since $(\partial h_{\osi}/\partial z
)({\mathbb H}^2)$ is a simply connected
domain in ${\mathbb C}$ not containing~$0$.).
With $\pi =\pi_1(M)$ and $\rho_{\mathbf w}$
as above, define a representation
\[ \rho^{\mathbb C}=\rho^{\mathbb C}_{\mathbf w}\co\pi\times {\mathbb Z}
\longrightarrow \widetilde{\SL}_2\times \E^1\]
by
\[ \rho^{\mathbb C}(u,0)=\rho_{\mathbf w} (u)\in\widetilde{\SL}_2\subset
\widetilde{\SL}_2\times \E^1\]
and, with $e_0$ denoting the neutral element of~$\pi$,
\[ \rho^{\mathbb C}(e_0,1)(z,w)=(z,w+1).\]
Since $h_{\sigma}$ commutes with translation in the $w$-coordinate and
the $\rho_{\mathbf w}(u)$ are lifts to $\widetilde{\ST {\mathbb H}^2}=
{\mathbb H}^2\times{\mathbb C}$ of
the differentials of the $\orh (\ou )$ (and hence unique up to
such translations), we can define a holomorphic map
\[ \sigma^{\mathbb C}(u,n)=\sigma_{\mathbf w}^{\mathbb C}(u,n)\co
D(\osi )\times {\mathbb C}\longrightarrow D(\osi )\times {\mathbb C}\]
for each $(u,n)\in\pi\times{\mathbb Z}$ such that $\sigma^{\mathbb C}(u,n)$
covers $\osi (\ou )$ and such that we get the commutative diagram
\[
\begin{CD}
{\mathbb H}^2\times {\mathbb C} @>\rho^{\mathbb C} (u,n)>>
       {\mathbb H}^2\times {\mathbb C}\\
@Vh_{\sigma}VV             @VVh_{\sigma}V\\
D(\osi )\times {\mathbb C} @>\sigma^{\mathbb C}(u,n)>>
                D(\osi )\times{\mathbb C}.
\end{CD}
\]
Notice that the choice of logarithm in the definition of $h_{\sigma}$ does
not affect the definition of $\sigma^{\mathbb C}(u,n)$, and that the map
\[ (\sigma ,z,w)\longmapsto (\sigma ,\sigma^{\mathbb C}(u,n)(z,w)),\;\;
\sigma\in\T (M),\; 
(z,w)\in D(\osi )\times {\mathbb C}\]
is holomorphic for any $(u,n)\in\pi\times{\mathbb Z}$.

Analogous to the construction in {\bf (1)} we define ${\mathcal E}_{\mathbb C}
(M)$ as the quotient space of
\[ \left\{ (\sigma ,z,w)\in \T (M)\times {\mathbb C}\times{\mathbb C}\co
z\in D(\osi )\right\} \]
under the equivalence relation
\begin{eqnarray*}
\lefteqn{(\sigma_1,z_1,w_1)\sim (\sigma_2,z_2,w_2)
:\Longleftrightarrow}\\
  &  & \sigma_1=\sigma_2\;\mbox{\rm and}\;
        (z_2,w_2)\in\sigma^{\mathbb C}(\pi\times{\mathbb Z})(z_1,w_1).
\end{eqnarray*}
With the complex structure thus induced on ${\mathcal E}_{\mathbb C}(M)$,
the projection
\[ \Psi_{\mathbb C}\co {\mathcal E}_{\mathbb C}(M)\longrightarrow \T (M)\]
is holomorphic and the fibre over $\sigma$ is the complex surface associated
(in the sense of~\cite{gego95}) to the taut contact circle on $M$ defined
by~$\sigma$. The universal family
\[ \Psi\co {\mathcal E}(M)\longrightarrow \T (M)\]
of taut contact circles on $M$ is then also defined as in~{\bf (1)}.

The action of $\Out (\pi )$ on $\T (M)$ extends to an action on
${\mathcal E}(M)$ such that ${\mathcal E}(M)/\Out (\pi )$ is a fibre space over
$\M (M)$. Here the fibre over $[\sigma ]\in\M (M)$ is the quotient of $M$
under the symmetries of the taut contact circle determined by~$\sigma$,
which correspond to the biholomorphisms of the complex surface $M\times S^1$
(defined by $\sigma^{\mathbb C}$) that are identical in the $S^1$-factor.

\vspace{2mm}

{\bf (3)} We end this paper with an alternative description of the moduli
of taut contact circles on~$S^3$. The corresponding statements for
left-quotients of $\SU (2)$ can be derived without much difficulty.

In the sequel, $(\omega_1,\omega_2)$ will always denote a pair of
pointwise linearly independent $1$-forms on an orientable $3$-manifold~$M$.
Then the line field $\ker\omega_1\cap\ker\omega_2$ is orientable
and thus admits nowhere zero sections.

\begin{defn}
{\rm (i)} The complex-valued $1$-form $\omega =\omega_1+\omega_2$ is
{\bf formally integrable} if $\omega\wedge d\omega\equiv 0$.

{\rm (ii)} A nowhere zero vector field $Y\in\ker\omega_1\cap\ker\omega_2$
defines a {\bf transversely conformal flow} if there is a complex-valued
function $f+ig$ on $M$ such that the Lie derivative of $\omega$ satisfies
$L_Y\omega =(f+ig)\omega$.
\end{defn}

Observe: (i) If $(\omega_1,\omega_2)$ is a taut contact circle, then
$\omega$ is formally integrable. (ii) The condition $L_Y\omega =
(f+ig)\omega$ implies that the flow of $Y$ preserves the conformal
structure transverse to $Y$ defined by $\omega_1\otimes\omega_1+
\omega_2\otimes\omega_2$.

The two definitions just given are in fact equivalent:

\begin{prop}
The complex $1$-form $\omega =\omega_1+i\omega_2$ is formally integrable if
and only if any nowhere zero vector field $Y\in\ker\omega_1
\cap\ker\omega_2$
defines a transversely conformal flow.
\end{prop}

This is immediate from the following equivalences:
\begin{eqnarray*}
\omega\wedge d\omega\equiv 0 & \Leftrightarrow & i_Y(\omega\wedge d\omega )
    \equiv 0 \;\;\mbox{\rm for some/every nonzero}\; Y\in\ker
    \omega_1\cap\ker\omega_2\\
 & \Leftrightarrow & \omega\wedge L_Y\omega\equiv 0\\
 & \Leftrightarrow & L_Y\omega =(f+ig)\omega.
\end{eqnarray*}

\vspace{2mm}

The proof of the following proposition is completely analogous to the
construction of the Godbillon-Vey invariant for real codimension~1
foliations, cf.~\cite{tamu92}, and will be left to the reader.

\begin{prop}
Assume $\omega =\omega_1+i\omega_2$ is formally integrable. Then there
is a complex $1$-form $\gamma$ such that $d\omega =\gamma\wedge\omega$,
and the complex number
\[ \int_M\gamma\wedge d\gamma \in{\mathbb C}\]
depends only on the oriented line field $\ker\omega_1\cap\ker\omega_2$
and the transverse conformal structure $\omega_1\otimes\omega_1+
\omega_2\otimes\omega_2$. In particular, if $(\omega_1,\omega_2)$ is a taut
contact circle, then this complex number is an invariant of the homothety
class of~$(\omega_1,\omega_2)$.
\end{prop}

\begin{defn}
We call the complex number $\int_M\gamma\wedge d\gamma$
the {\bf Godbillon-Vey invariant}
of the transversely conformal flow defined by~$\omega$.
\end{defn}

Write $|S^3|$ for the volume of $S^3$ with repect to the volume form
\[ \sum_{i=1}^4 x_i\partial_{x_i}\rfloor (dx_1\wedge dx_2
\wedge dx_3\wedge dx_4). \]
For the taut contact circle corresponding to $[a]\in \M_1$ in the notation
of Theorem~\ref{thm:deform-s}, i.e.\ the class of $a$ and $1-a$, the
Godbillon-Vey invariant can be computed as
\[ GV([a])=\int_{S^3}\gamma_a\wedge d\gamma_a =-\frac{2|S^3|}{a(1-a)}.\]
For the taut contact circle corresponding to $n\in\M_2$ one finds
\[ GV(n)=\int_{S^3}\gamma_n\wedge d\gamma_n =-\frac{2n|S^3|}{(n+1)^2}.\]
This will be proved in~\cite{gegoc}.

Notice that the class $[a]\in\M_1$ is determined by the complex number
$a(1-a)$. So the Godbillon-Vey invariant defines a biholomorphism
\[ \begin{array} {rcl}
\M_1 & \longrightarrow & \{ (x+iy\in{\mathbb C}\co x\geq y^2\} \\
\mbox{$[a]$} & \longmapsto & a(1-a)=-2|S^3|/GV([a]).
\end{array} \]

Other aspects of the Godbillon-Vey invariant will be discussed in~\cite{gegoc}.
For instance, for the pair of Liouville-Cartan forms on the unit
cotangent bundle of a surface with a Riemannian metric, the Godbillon-Vey
invariant equals the total Gau{\ss} curvature of the metric, and one can prove
a generalisation of the Gau{\ss}-Bonnet theorem about the topological
invariance of this total curvature, showing that under suitable assumptions
the Godbillon-Vey invariant is not just an invariant of the transversely
conformal structure, but in fact only depends on the oriented and
cooriented common kernel foliation.

Furthermore, we shall prove there that the homothety class of a taut contact
circle on $S^3$ is determined, up to diffeomorphism, by the oriented
and cooriented common kernel foliation, and define a diffeomorphism invariant
of such foliations that allow directly to recover the modulus of the taut
contact circle.

\vspace{2mm}

\noindent {\bf Acknowledgements.} This work was initiated during several brief
visits by J.~G. to the Forschungsinstitut f\"ur Mathematik
(FIM) of the ETH Z\"urich between 1995 and 1997, while H.~G. held an assistant
professorship there, and a visit by  H.~G. to Madrid, supported by
DGICYT (Direcci\'on General de Investigaci\'on Cient\'{\i}fica y
T\'ecnica)
grant no.\ PB95-0185. A major share of the research was done during a visit
by J.~G. to Leiden University from March till May 2000, supported by
an NWO (Nederlandse Organisatie voor Wetenschappelijk Onderzoek)
bezoekersbeurs. The final version was prepared while H.~G. was a guest
of the Universidad Aut\'onoma.
We thank these funding organisations and universities for their support.

We also thank Horst Kn\"orrer for helpful discussions.



\noindent Hansj\"org Geiges, Mathematisch Instituut, Universiteit Leiden,
Postbus 9512, 2300 RA Leiden, The Netherlands; e-mail:
{\tt geiges@math.leidenuniv.nl}

\vspace{1mm}

\noindent Jes\'us Gonzalo, Departamento de Matem\'aticas, Universidad
Aut\'onoma de Madrid, 28049 Madrid, Spain; e-mail:
{\tt jesus.gonzalo@uam.es}
\end{document}